# On finite homomorphic images of the multiplicative group of a division algebra

By Yoav Segev*

### Introduction

The purpose of this paper, together with [6], is to prove that the following Conjecture 1 holds:

CONJECTURE 1 (A. Potapchik and A. Rapinchuk). *Let $D$ be a finite dimensional division algebra over an arbitrary field. Then $D^{\#}$ does not have any normal subgroup $N$ such that $D^{\#}/N$ is a nonabelian finite simple group.*

Of course $D^{\#}$ is the multiplicative group of $D$. Conjecture 1 appears in [4]. It is related to the following conjecture of G. Margulis and V. Platonov (Conjectures 9.1 and 9.2, pages 510–511 in [3], or Conjecture (PM) in [4]).

CONJECTURE 2 (G. Margulis and V. Platonov). *Let $\mathfrak{G}$ be a simple, simply connected algebraic group defined over an algebraic number field $K$. Let $T$ be the set of all nonarchimedean places $v$ of $K$ such that $\mathfrak{G}$ is $K_v$-anisotropic; then for any noncentral normal subgroup $N \leq \mathfrak{G}(K)$ there exists an open normal subgroup $W \leq \mathfrak{G}(K,T) = \prod_{v \in T} \mathfrak{G}(K_v)$ such that $N = \mathfrak{G}(K) \cap W$; in particular, if $T = \emptyset$ then $\mathfrak{G}(K)$ does not have proper noncentral normal subgroups.*

In Corollary 2.5 of [4], Potapchik and Rapinchuk prove that if $D$ is a finite dimensional division algebra over an algebraic number field $K$, then for $\mathfrak{G} = \mathrm{SL}_{1,D}$, Conjecture 2 is equivalent to the nonexistence of a normal subgroup $N \triangleleft D^{\#}$ such that $D^{\#}/N$ is a nonabelian finite simple group. Of course this was the main motivation for the conjecture of Potapchik and Rapinchuk in [4]. Thus as a corollary, we get that if $D$ is a finite dimensional division algebra over an algebraic number field $K$ and $\mathfrak{G} = \mathrm{SL}_{1,D}$, then the normal subgroup structure of $\mathfrak{G}(K)$ is given by Conjecture 2.

Hence we prove Conjecture 2, in one of the cases when $\mathfrak{G}$ is of type $A_n$. The case when $\mathfrak{G}$ is of type $A_n$ is the main case left open in Conjecture 2. For

*This work was partially supported by BSF 92-003200 and by grant no. 6782-1-95 from the Israeli Ministry of Science and Art.

further information about the historical background and the current state of Conjecture 2, we refer the reader to Chapter 9 in [3] and to the introduction in [4].

More generally we are interested in the possible structure of finite homomorphic images of the multiplicative group of a division algebra. Let $D$ be a division algebra and let $D^\#$ denote the multiplicative group of $D$. Various papers dealt with subgroups of finite index in $D^\#$, e.g., [2], [4], [7] and the references therein. We refer the reader to [1], for a survey article on the history of finite dimensional central division algebras.

Let $X$ be a finite group. Define the commuting graph of $X$, $\Delta(X)$ as follows. Its vertex set is $X \setminus \{1\}$. Its edges are pairs $\{a, b\}$, such that $a, b \in X \setminus \{1\}$, $a \neq b$, and $[a, b] = 1$ ($a$ and $b$ commute). We denote the diameter of $\Delta(X)$ by $\mathrm{diam}(\Delta(X))$.

Let $\mathrm{d} : \Delta(X) \times \Delta(X) \to \mathbb{Z}^{\geq 0}$ be the distance function on $\Delta(X)$. We say that $\Delta(X)$ is balanced if there exist $x, y \in \Delta(X)$ such that the distances $d(x, y), d(x, xy), d(y, xy), d(x, x^{-1}y), d(y, x^{-1}y)$ are all bigger than 3.

The Main Theorem of this paper is:

THEOREM A. *Let $L$ be a nonabelian finite simple group. Suppose that either $\mathrm{diam}(\Delta(L)) > 4$, or $\Delta(L)$ is balanced. Let $D$ be a finite dimensional division algebra over an arbitrary field. Then $D^\#$ does not have any normal subgroup $N$ such that $D^\#/N \simeq L$.*

The proof of Theorem A does not rely on the classification of finite simple groups. However, in [6] we prove (using classification) that all nonabelian finite simple groups $L$ have the property that $\Delta(L)$ is balanced or $\mathrm{diam}(\Delta(L)) > 4$. Thus Theorem A together with [6] prove the assertion of Conjecture 1.

The *organization of the proof of Theorem* A is as follows. Let $D$ be a division algebra (not necessarily finite dimensional over its center $F := Z(D)$). Let $G := D^\#$ be the multiplicative group of $D$ and let $N$ be a normal subgroup of $G$ such that $G^* := G/N$ is finite (not necessarily simple). Let $\Delta = \Delta(G^*)$ be the commuting graph of $G^*$.

In Section 1 we introduce some notation and preliminaries. In particular we introduce the set $N(a)$, for $a \in G$, which plays a crucial role in the paper. In Section 2 we deal with $\Delta$ and note that severe restrictions are imposed on $\Delta$.

In Section 3 we introduce the *U-Hypothesis* which plays a central role throughout the paper. In addition, we establish in Section 3 some notation and preliminary results regarding the $U$-Hypothesis and we prove that if $\mathrm{diam}(\Delta) > 4$, then $G$ satisfies the $U$-Hypothesis. In Section 4 we show that if $\Delta$ is balanced then $G$ satisfies the $U$-Hypothesis. Sections 5 and 6 are independent of the rest of the paper and deal with further consequences of the $U$-Hypothesis.

From Section 7 to the end of the paper, we specialize to the case when $D$ is finite dimensional over $F$ and $G^*$ is nonabelian simple. We assume that either $\operatorname{diam}(\Delta) > 4$, or $\Delta$ is balanced and set out to obtain our contradiction. Section 7 gives some preliminaries and technical results. In particular, we introduce in Section 7 (see the definitions at the beginning) *the set $\hat{K}$, which plays a crucial role in the proof.* Sections 8 and 9 are basically devoted to the proof that $\hat{K} = \mathbb{O}U \setminus N$ (Theorem 9.1), which is the main target of the paper. Once Theorem 9.1 is proved, we can use it in Section 10 to construct a local ring $R$, whose existence yields a contradiction and proves Theorem A.

## 1. Notation and preliminaries

All through this paper $D$ is a division algebra over its center $F := Z(D)$. In some sections we will assume that $D$ is finite dimensional over $F$, but in general we do not assume this. We let $D^\# = D \setminus \{0\}$ and $G = D^\#$ be the multiplicative group of $D$. Letting $F^\# = F \setminus \{0\}$, we denote $N$ a normal subgroup of $G$ such that $F^\# \leq N$ and $G/N$ is finite. The following notational convention is used: $G^* = G/N$ and for $a \in G$, we let $a^*$ denote its image in $G^*$ under the canonical homomorphism; that is, $a^* = Na$. If $H^*$ is a subgroup of $G^*$, then by convention $H \leq G$ is the full inverse image of $H^*$ in $G$.

(1.1) *Remark.* Note that since $F^\# \leq N$, for all $a \in G$ and $\alpha \in F^\#$, $(\alpha a)^* = a^*$, and in particular, $(-a)^* = a^*$. We use this fact without further reference.

(1.2) *Notation.* (1) Let $a \in G$. We denote
$$N(a) = \{n \in N : a + n \in N\}.$$
(2) Let $A, B \subseteq D$. We denote $A + B = \{a + b : a \in A, b \in B\}$, $A - B = \{a - b : a \in A, b \in B\}$ and $-A = \{-a : a \in A\}$.
(3) Let $A, B \subseteq D$ and $x \in D$. We denote $AB = \{ab : a \in A, b \in B\}$, $Ax = \{ax : a \in A\}$ and $xA = \{xa : a \in A\}$.
(4) We denote by $[D : F]$ the dimension of $D$ as a vector space over $F$. If $[D : F] < \infty$, then as is well known $[D : F] = n^2$, for some natural $n \geq 1$. We denote $\deg(D) = n$.

(1.3) *Notation for the case $[D : F] < \infty$.* If $[D : F] < \infty$, we denote

(1) $$\nu : G \to F^\#$$

the reduced-norm function. Of course $\nu$ is a group homomorphism.

(2) $$\mathbb{O} = \mathbb{O}(D) = \{a \in D^\# : \nu(a) = 1\}.$$

(1.4) *Suppose $[D : F] < \infty$. Then for all $a \in G$, $\nu(a)$ is a product of conjugates of $a$ in $G$.*

*Proof.* This is well known and follows from Wedderburn's Factorization Theorem. See, e.g., [5, p. 253].

(1.5) *If $[D : F] < \infty$ and $[G^*, G^*] = G^*$, then $G = \mathbb{O}N$.*

*Proof.* Since $G/\mathbb{O}$ is isomorphic to a subgroup of $F^{\#}$, $G/\mathbb{O}$ is abelian, and hence $G/\mathbb{O}N$ is abelian. But $G/\mathbb{O}N \simeq (G/N)/(\mathbb{O}N/N)$, and hence $G^* = [G^*, G^*] \leq \mathbb{O}N/N$. Hence $G = \mathbb{O}N$.

(1.6) THEOREM (G. Turnwald). *Let $\mathfrak{D}$ be an infinite division algebra. Let $H \leq \mathfrak{D}^{\#}$ be a subgroup of finite index. Then $\mathfrak{D} = H - H$.*

*Proof.* This is a special case of Theorem 1 in [7].

(1.7) COROLLARY. $N + N = D = N - N$.

*Proof.* This follows from 1.6. Note that as $-1 \in N$, $N + N = N - N$.

(1.8) *Let $a \in G \setminus N$ and let $n \in N$. Then*
(1) $N(na) = nN(a)$ and $N(an) = N(a)n$.
(2) *For all $b \in G$, $N(b^{-1}ab) = b^{-1}N(a)b$.*
(3) $N(a) \neq \emptyset$.
(4) *If $n \in N(a)$, then $n^{-1} \notin N(a^{-1})$.*
(5) *There exists $a' \in Na$, with $1 \in N(a')$.*

*Proof.* In (1), we prove that $N(na) = nN(a)$. The proof that $N(an) = N(a)n$ is similar. Let $m \in N(na)$. Then $na + m \in N$. Hence $a + n^{-1}m \in N$, so $n^{-1}m \in N(a)$. Hence $m \in nN(a)$. Let $m \in nN(a)$. Then there exists $s \in N(a)$ such that $m = ns$. Then $na + m = na + ns = n(a + s)$. Since $s \in N(a)$, $a + s \in N$, so $na + m \in N$. Hence $m \in N(na)$.

For (2), let $m \in N(b^{-1}ab)$. Then $b^{-1}ab + m \in N$, and hence $a + bmb^{-1} \in N$. Hence $bmb^{-1} \in N(a)$, so $m \in b^{-1}N(a)b$. Let $m \in b^{-1}N(a)b$. Then there exists $s \in N(a)$, with $m = b^{-1}sb$. Then $b^{-1}ab + m = b^{-1}ab + b^{-1}sb = b^{-1}(a+s)b \in N$. Thus $m \in N(b^{-1}ab)$.

For (3), note that by 1.7 there exists $m, n \in N$ such that $a = n - m$. Hence $m \in N(a)$. Let $n \in N(a)$. Then $a + n \in N$. Multiplying by $a^{-1}$ on the right and by $n^{-1}$ on the left we get that $a^{-1} + n^{-1} \in Na^{-1}$, hence $n^{-1} \notin N(a^{-1})$. This proves (4). Finally to prove (5), let $n \in N(a)$. Then $1 \in n^{-1}N(a) = N(n^{-1}a)$.

(1.9) *Let $K$ be a finite group and let $\emptyset \neq \mathcal{A} \subsetneq K$ be a proper normal subset of $K$. Set $X := \{x \in K : x\mathcal{A} \subseteq \mathcal{A}\}$. Then $X$ is a proper normal subgroup of $K$. In particular, if $X \neq 1$, then $K$ is not simple.*

*Proof.* Since $\mathcal{A}$ is finite, $X = \{x \in K : x\mathcal{A} = \mathcal{A}\}$. Hence clearly $X$ is a subgroup of $K$. Let $y \in K$ and $x \in X$; then $(y^{-1}xy)\mathcal{A} = (y^{-1}xy)(y^{-1}\mathcal{A}y) = y^{-1}(x\mathcal{A})y = y^{-1}\mathcal{A}y = \mathcal{A}$, since $\mathcal{A}$ is a normal subset of $K$. Hence $y^{-1}xy \in X$, so $X$ is a normal subgroup of $K$. Clearly since $\mathcal{A}$ is a proper nonempty subset, $X \neq G$.

## 2. The commuting graph of $G^*$

Throughout the paper we let $\Delta$ be the graph whose vertex set is $G^* \setminus \{1^*\}$ and whose edges are $\{a^*, b^*\}$ such that $[a^*, b^*] = 1^*$. We call $\Delta$ the commuting graph of $G^*$ and let $\mathrm{d} : \Delta \times \Delta \to \mathbb{Z}^{\geq 0}$ be the distance function of $\Delta$.

(2.1) *Let $a \in G \setminus N$ and $n \in N$. Suppose that $a + n \in G \setminus N$. Let $H \leq G$, with $H^* = C_{G^*}(a^*)$. Then $(a+n)^* \in H^*$, so $a + n \in H$.*

*Proof.* Note that $n^{-1}a + 1 \in C_G(n^{-1}a)$. Thus $(n^{-1}a + 1)^* \in C_{G^*}((n^{-1}a)^*)$ $= C_{G^*}(a^*)$. But since $a + n = n(n^{-1}a + 1)$, $(a+n)^* = (n^{-1}a + 1)^*$.

(2.2) *Remark.* Note that by 2.1, if $a, b \in G \setminus N$ and $n \in N$, then if $a + b \in N$, or $a - b \in N$, $\mathrm{d}(a^*, b^*) \leq 1$ and if $n \notin N(a)$, then $\mathrm{d}((a+n)^*, a^*) \leq 1$. We use these facts without further reference.

(2.3) *Let $a, b, c \in G \setminus N$, with $a + b = c$. Then*
(1) *If $\mathrm{d}(a^*, b^*) > 2$, then $N(c) \subseteq N(a) \cap N(b)$.*
(2) *If $\mathrm{d}(a^*, b^*) > 2$, and $\mathrm{d}(a^*, c^*) > 2$, then $N(b) = N(c) \subseteq N(a) \cap N(-a)$.*
(3) *If $\mathrm{d}(a^*, b^*) > 4$, then either $N(a) = N(c) \subseteq N(b) \cap N(-b)$, or $N(b) = N(c) \subseteq N(a) \cap N(-a)$.*

*Proof.* For (1), let $n \in N(c) \setminus (N(a) \cap N(b))$. Suppose $n \notin N(a)$. Then
$$c + n = (a + n) + b.$$
As $c + n \in N$, 2.2 implies that $\mathrm{d}(a^*, (a+n)^*) \leq 1 \geq \mathrm{d}(b^*, (a+n)^*)$; thus $\mathrm{d}(a^*, b^*) \leq 2$, a contradiction.

Assume the hypotheses of (2). By (1), $N(c) \subseteq N(a) \cap N(b)$ and since $b = c - a$, (1) implies that $N(b) \subseteq N(c) \cap N(-a)$. Hence (2) follows. (3) follows from (2) since we must have either $\mathrm{d}(a^*, c^*) > 2$, or $\mathrm{d}(b^*, c^*) > 2$.

(2.4) *Remark.* Note that by 2.3.3, if $a, b \in G \setminus N$, with $\mathrm{d}(a^*, b^*) > 4$, then $N(a) \subseteq N(b)$, or $N(b) \subseteq N(a)$. We use this fact without further reference.

(2.5) *Let $a, b \in G \setminus N$ such that $d(a^*, b^*) > 1$ and $N(a) \not\subseteq N(b)$. Then*
(1) $b^*(b+n)^*(a-b)^*$ *is a path in $\Delta$, for any $n \in N(a) \setminus N(b)$.*
(2) *If $-1 \notin N(ab^{-1})$, then for all $n \in N(a) \setminus N(b)$,*

$$b^*(b+n)^*(ab^{-1}-1)^*(ab^{-1})^* \text{ is a path in } \Delta.$$

(3) *If $-1 \notin N(b^{-1}a)$, then for all $n \in N(a) \setminus N(b)$,*

$$b^*(b+n)^*(b^{-1}a-1)^*(b^{-1}a)^* \text{ is a path in } \Delta.$$

*Proof.* Let $c = a - b$. Since $d(a^*, b^*) > 1$, $c \notin N$. Next note that $c + b = a$. Let $n \in N(a) \setminus N(b)$. Then $c + (b+n) = a + n \in N$. Hence $d(c^*, (b+n)^*) \leq 1$. This show (1).

Suppose $-1 \notin N(ab^{-1})$ and let $n \in N(a) \setminus N(b)$. Note that $c = (ab^{-1}-1)b$. Further, $c^*$ commutes with $(b+n)^*$ and $b^*$ commutes with $(b+n)^*$. It follows that $d((ab^{-1}-1)^*, (b+n)^*) \leq 1$. Clearly $d((ab^{-1}-1)^*, (ab^{-1})^*) \leq 1$, so (2) follows. The proof of (3) is similar to the proof of (2) when we notice that $c = b(b^{-1}a - 1)$.

(2.6) *Let $a, b \in G \setminus N$ with $d(a^*, b^*) > 4$. Suppose $N(a) \subseteq N(b)$. Then*
(1) $N(a+b) = N(a) \subseteq N(b) \cap N(-b)$.
(2) $N(a-b) = N(a) \subseteq N(b) \cap N(-b)$.

*Proof.* For (1) we use 2.3.3. Suppose (1) is false. Set $c = a + b$. Then by 2.3.3, $N(b) = N(c) \subseteq N(a) \cap N(-a)$. Since $N(a) \subseteq N(b)$, we must have $N(b) = N(c) = N(a) \cap N(-a) = N(a)$. It follows that $N(a) \subseteq N(-a) = -N(a)$. Multiplying by $-1$, we get that $N(-a) \subseteq N(a)$, so $N(a) = N(-a)$. Thus $N(b) = N(c) = N(a) = N(-a)$. Hence $N(a) = N(c) \subseteq N(b) \cap N(-b)$ in this case too.

Suppose (2) is false. Set $c = a - b$. Then by 2.3.3, $N(-b) = N(c) \subseteq N(a) \cap N(-a)$. In particular $N(-b) \subseteq N(-a)$, so $N(b) \subseteq N(a)$. Hence we must have $N(-b) = N(c) = N(a) \cap N(-a) = N(-a)$. As above we get that $N(a) = N(-a) = N(b) = N(c)$, so again $N(c) = N(a) \subseteq N(b) \cap N(-b)$.

(2.7) *Let $a, b \in G \setminus N$. Suppose*
(a) $d(a^*, b^*) > 4$.
(b) $N(a) \subseteq N(b)$.
*Then*
(1) *If $1 \in N(a)$, then $\pm 1 \in N(b)$.*
(2) *For all $n \in N \setminus N(b)$*

$$N(a) \subseteq N(a+n) \text{ and } -N(a) \subseteq N(b+n) \supseteq N(a).$$

*Proof.* Set $x = a - b$. Note first that by 2.6.2,

$$(*) \qquad N(a) = N(x) \subseteq N(b) \cap N(-b).$$

Note that this already implies (1). Next note that

$$x = (a + n) - (b + n).$$

Since $d(a^*, b^*) > 4$, we get that $d((a+n)^*, (b+n)^*) > 2$. Hence by 2.3.1, $N(x) \subseteq N(a+n) \cap N(-(b+n))$. Thus $N(a) = N(x) \subseteq N(a+n)$ and $N(a) = N(x) \subseteq N(-(b+n))$, so that $-N(a) \subseteq N(b+n)$.

Finally, note that by $(*)$, $N(-a) \subseteq N(b)$, so by the previous paragraph of the proof $-N(-a) \subseteq N(b+n)$, that is $N(a) \subseteq N(b+n)$ and the proof of 2.7 is complete.

(2.8) *Let $a, b \in G \setminus N$ be such that $ab \in G \setminus N$. Then*
(1) *Assume $N(ab) \not\supseteq N(b)$ and $-1 \notin N(a^{-1})$. Then for all $m \in N(b) \setminus N(ab)$,*

$$a^*(a^{-1} - 1)^*(ab + m)^*(ab)^* \text{ is a path in } \Delta.$$

(2) *Assume $N(ab) \not\supseteq N(a)$, and $-1 \notin N(b^{-1})$; then for all $m \in N(a) \setminus N(ab)$*

$$b^*(b^{-1} - 1)^*(ab + m)^*(ab)^* \text{ is a path in } \Delta.$$

*Proof.* We have

$$(1 - a)b + ab = b.$$

Let $m \in N(b) \setminus N(ab)$. Then

$$(1 - a)b + ab + m = b + m \in N.$$

This implies that $(ab+m)^*$ commutes with $(1-a)^*b^*$. Of course $(ab+m)^*$ commutes also with $a^*b^*$. Hence $(ab+m)^*$ commutes with $((1-a)^*b^*)(b^*)^{-1}(a^*)^{-1}$ $= (a^{-1}-1)^*$. Hence we conclude that $a^*(a^{-1}-1)^*(ab+m)^*(ab)^*$ is a path in $\Delta$, this completes the proof of (1). The proof of (2) is similar since $a(1-b)+ab = a$.

(2.9) *Let $a, b \in G \setminus N$. Then*
(1) *Assume that $N(ab) \not\subseteq N(a)$ and $-1 \notin N(b)$. Then for all $m \in N(ab) \setminus N(a)$,*

$$a^*(a + m)^*(b - 1)^*b^* \text{ is a path in } \Delta.$$

(2) *Assume that $N(ab) \not\subseteq N(b)$ and $-1 \notin N(a)$. Then for all $m \in N(ab) \setminus N(b)$,*

$$a^*(a - 1)^*(b + m)^*b^* \text{ is a path in } \Delta.$$

*Proof.* First note that
$$a(b-1) + a = ab.$$
Let $m \in N(ab) \setminus N(a)$. Then
$$a(b-1) + a + m = ab + m \in N.$$
Hence $(a+m)^*$ commutes with $a^*(b-1)^*$. Of course $(a+m)^*$ commutes with $a^*$, so $(a+m)^*$ commutes with $(b-1)^*$. Hence $a^*(a+m)^*(b-1)^*b^*$ is a path in $\Delta$. This proves (1). The proof of (2) is similar because $(a-1)b + b = ab$.

(2.10) *Let* $a, b \in G \setminus N$. *Assume*
(i) $-1 \notin N(a) \cup N(b)$.
(ii) *For all* $g \in G$, $-1 \in N(ab^g)$.
*Then* $G^*$ *is not simple.*

*Proof.* Let $g \in G$. Note that by 1.8.2, $-1 \notin N(b^g)$, for all $g \in G$. Thus by (ii), $N(ab^g) \not\subseteq N(b^g)$ and $-1 \in N(ab^g) \setminus N(b^g)$. Hence by 2.9.2, $a^*(a-1)^*(b^g-1)^*b^*$ is a path in $\Delta$. In particular

(∗) $\qquad\qquad d((a-1)^*, (b^g-1)^*) \leq 1$, for all $g \in G$.

Note now that $(b^g-1)^* = ((b-1)^*)^{g^*}$, so that $C^* := \{(b^g-1)^* : g \in G\}$ is a conjugacy class of $G^*$. Now (∗) implies that $(a-1)^*$ commutes with every element of $C^*$, so that $G^*$ is not simple.

(2.11) *Let* $x, y \in G \setminus N$ *and* $n, m \in N$ *such that*
(a) $xny \notin N$.
(b) $m \in N(xn) \cap N(ny)$.
(c) $-1 \notin N(ny) \cap N(xn)$.
(d) $m \notin N(x) \cup N(y)$.
(e) $-1 \notin N(x^{-1}) \cup N(y^{-1})$.
*Then*
(1) *If* $m \in N(xny)$ *and* $-1 \notin N(ny)$, *then* $x^*(x+m)^*(ny-1)^*y^*$ *is a path in* $\Delta$.
(2) *If* $m \in N(xny)$ *and* $-1 \notin N(xn)$, *then* $x^*(xn-1)^*(y+m)^*y^*$ *is a path in* $\Delta$.
(3) *If* $m \notin N(xny)$, *then* $x^*(x^{-1}-1)^*(xny+m)^*(y^{-1}-1)^*y^*$ *is a path in* $\Delta$.
(4) $d(x^*, y^*) \leq 4$.

*Proof.* Suppose first that $m \in N(xny)$ and $-1 \notin N(ny)$; then since $m \notin N(x)$, we see that $m \in N(xny) \setminus N(x)$. Since $-1 \notin N(ny)$, we get (1) from 2.9.1.

Suppose next that $m \in N(xny)$ and $-1 \notin N(xn)$; then since $m \notin N(y)$, we see that $m \in N(xny) \setminus N(y)$. Since $-1 \notin N(xn)$, we get (2) from 2.9.2.

Now assume $m \notin N(xny)$. Since $m \in N(xn)$, we see that $m \in N(xn) \setminus N(xny)$. Further, $-1 \notin N(y^{-1})$; hence, by 2.8.2, $y^*(y^{-1} - 1)^*(xny + m)^*$ is a path in $\Delta$. Next, since $m \in N(ny)$, we see that $m \in N(ny) \setminus N(xny)$. Further $-1 \notin N(x^{-1})$; hence, by 2.8.1, $x^*(x^{-1} - 1)^*(xny + m)^*$ is a path in $\Delta$. Hence (3) follows and (4) is immediate from (1), (2) and (3).

## 3. The definition of the $U$-Hypothesis; notation and preliminaries; the proof that if $\mathrm{diam}(\Delta) > 4$ then $G$ satisfies the $U$-Hypothesis

In this section we define the *U-Hypothesis* which will play a crucial role in the paper. We also establish some notation which will hold throughout the paper and give some preliminary results. Finally, in Theorem 3.18, we prove that if $\mathrm{diam}(\Delta) > 4$, then $G$ satisfies the $U$-Hypothesis.

*Definition.* We say that $G$ satisfies the *U-Hypothesis* with respect to $\mathbb{N}$ (or just that $G$ satisfies the $U$-Hypothesis) if there exists a normal subset $\emptyset \neq \mathbb{N} \subsetneq G$ such that $\mathbb{N} \subsetneq N$ is a proper subset of $N$ and if we set $\bar{\mathbb{N}} = N \setminus \mathbb{N}$, then

(U1) $1, -1 \in \mathbb{N}$.
(U2) $\mathbb{N}^2 = \mathbb{N}$.
(U3) For all $\bar{n} \in \bar{\mathbb{N}}$, $\bar{n} + 1 \in \mathbb{N}$ and $\bar{n} - 1 \in N$.

*Notation.* Let $x^* \in G^* \setminus \{1^*\}$ and let $C^* \subseteq G^* - \{1^*\}$ be a conjugacy class of $G^*$.

(1) Denote $\mathbb{P}_{x^*} = \{a \in Nx : 1 \in N(a)\}$.
(2) Denote
$$\mathbb{N}_{x^*} = \{n \in N : n \in N(a), \text{ for all } a \in \mathbb{P}_{x^*}\},$$
$$\bar{\mathbb{N}}_{x^*} = N \setminus \mathbb{N}_{x^*}.$$

(3) Let $U_{x^*} = \{n \in N : n, n^{-1} \in \mathbb{N}_{x^*}\}$.
(4) Let $\mathbb{M}_{x^*} = \mathbb{N}_{x^*} \setminus U_{x^*}$.
(5) Let $\mathbb{O}_{x^*} = \{x_1 \in Nx : -1 \notin N(x_1) \cup N(x_1^{-1})\}$.
(6) Denote by $C_{x^*}$ the conjugacy class of $x^*$ in $G^*$.
(7) Denote $\hat{C} = \{c \in G : c^* \in C^*\}$.
(8) Let $\mathbb{P}_{C^*} = \bigcup_{y^* \in C^*} \mathbb{P}_{y^*}$.
(9) Denote
$$\mathbb{N}_{C^*} = \bigcap_{y^* \in C^*} \mathbb{N}_{y^*},$$
$$\bar{\mathbb{N}}_{C^*} = N \setminus \mathbb{N}_{C^*}.$$

(10) Denote $U_{C^*} = \bigcap_{y^* \in C^*} U_{y^*} = \{n \in N : n, n^{-1} \in \mathbb{N}_{C^*}\}$.

(11) Let $\mathbb{M}_{C^*} = \mathbb{N}_{C^*} \setminus U_{C^*}$.

*Definition.* We define three binary relations on $(G^* \setminus \{1^*\}) \times (G^* \setminus \{1^*\})$. These relations will play a crucial role throughout this paper. Given a binary relation $R$ on $(G^* \setminus \{1^*\}) \times (G^* \setminus \{1^*\})$, $R(x^*, y^*)$ means that $(x^*, y^*) \in R$. Here are our binary relations: Let $(x^*, y^*) \in (G^* \setminus \{1^*\}) \times (G^* \setminus \{1^*\})$.

$\operatorname{In}(x^*, y^*)$: For all $a \in Nx$ and $b \in Ny$, either $N(a) \subseteq N(b)$, or $N(b) \subseteq N(a)$. Note that $\operatorname{In}(x^*, y^*)$ is a symmetric relation.

$\operatorname{Inc}(y^*, x^*)$: $\operatorname{In}(y^*, x^*)$ and for all $b \in \mathbb{P}_{y^*}$, there exists $a \in \mathbb{P}_{x^*}$ such that $N(b) \supseteq N(a)$. Note that $\operatorname{Inc}(y^*, x^*)$ is not necessarily symmetric.

$T(x^*, y^*)$: For all $(a, b) \in Nx \times Ny$, and all $n \in N \setminus (N(a) \cup N(b))$
$$N(a+n) \supseteq N(a) \cap N(b) \subseteq N(b+n).$$

Note that $T(x^*, y^*)$ is symmetric.

(3.1) *Let $x^*, y^* \in G^* \setminus \{1^*\}$ and let $g \in G$. Then*

(1) $g^{-1}\mathbb{P}_{x^*}g = \mathbb{P}_{(g^{-1}xg)^*}$.
(2) $g^{-1}\mathbb{N}_{x^*}g = \mathbb{N}_{(g^{-1}xg)^*}$ *and* $g^{-1}\bar{\mathbb{N}}_{x^*}g = \bar{\mathbb{N}}_{(g^{-1}xg)^*}$.
(3) $\mathbb{N}_{C_{x^*}}$ *is a normal subset of $G$.*
(4) *If* $-1 \in \mathbb{N}_{x^*}$, *then* $-1 \in \mathbb{N}_{C_{x^*}}$.
(5) *If* $\mathbb{N}_{y^*} \supseteq \mathbb{N}_{x^*}$, *then* $\mathbb{N}_{C_{y^*}} \supseteq \mathbb{N}_{C_{x^*}}$.
(6) $g^{-1}\mathbb{M}_{x^*}g = \mathbb{M}_{(g^{-1}xg)^*}$, $g^{-1}U_{x^*}g = U_{(g^{-1}xg)^*}$ *and* $g^{-1}\mathbb{O}_{x^*}g = \mathbb{O}_{(g^{-1}xg)^*}$.

*Proof.* For (1), it suffices to show that $g^{-1}\mathbb{P}_{x^*}g \subseteq \mathbb{P}_{(g^{-1}xg)^*}$. Let $a \in \mathbb{P}_{x^*}$. Then $a \in Nx$ and $1 \in N(a)$, so that, by 1.8, $1 \in N(a^g)$, and clearly, $a^g \in Nx^g$. Hence $a^g \in \mathbb{P}_{(x^g)^*}$. For (2), it suffices to show that $g^{-1}\mathbb{N}_{x^*}g \subseteq \mathbb{N}_{(g^{-1}xg)^*}$. Let $n \in \mathbb{N}_{x^*}$. Then $n \in N(a)$, for all $a \in \mathbb{P}_{x^*}$; hence, by 1.8, $n^g \in N(c)$, for all $c \in g^{-1}\mathbb{P}_{x^*}g$. Now, by (1), $n^g \in \mathbb{N}_{(g^{-1}xg)^*}$. Note that (3) and (4) are immediate from (2).

For (5), let $z^* \in C_{y^*}$. Let $g \in G$, with $(y^g)^* = z^*$. By (2), $\mathbb{N}_{z^*} \supseteq \mathbb{N}_{(x^g)^*} \supseteq \mathbb{N}_{C_{x^*}}$. As this holds for all $z^* \in C_{y^*}$, we see that $\mathbb{N}_{C_{y^*}} \supseteq \mathbb{N}_{C_{x^*}}$.

The proof of (6) is similar to the proof of (2) and we omit the details.

(3.2) *Let $x^* \in G^* \setminus \{1^*\}$, let $\alpha \in \{x^*, C_{x^*}\}$ and set $\mathbb{P} = \mathbb{P}_\alpha$ and $\mathbb{N} = \mathbb{N}_\alpha$. Then*

(1) $1 \in \mathbb{N}$.
(2) $n \in \mathbb{N}$ *if and only if* $n^{-1}\mathbb{P} \subseteq \mathbb{P}$.
(3) *If* $n \in \mathbb{N}$, *then* $n\mathbb{N} \subseteq \mathbb{N}$.
(4) *If* $\alpha = C_{x^*}$, *then* $\mathbb{N}$ *is a normal subset of $G$.*
(5) *If* $-1 \in \mathbb{N}$, *then* $-\mathbb{N} = \mathbb{N}$.

*Proof.* (1) is by the definition of $\mathbb{N}$. Let $n \in N$. Suppose $n^{-1}\mathbb{P} \subseteq \mathbb{P}$. Let $a \in \mathbb{P}$. Then $n^{-1}a \in \mathbb{P}$ and hence, $1 \in N(n^{-1}a)$; so by 1.8.1, $n \in N(a)$. As this holds for all $a \in \mathbb{P}$, $n \in \mathbb{N}$. Suppose $n \in \mathbb{N}$ and let $a \in \mathbb{P}$; then $n \in N(a)$; so by 1.8.1, $1 \in N(n^{-1}a)$, and $n^{-1}a \in \mathbb{P}$.

Let $n \in \mathbb{N}$. Then by (2), for all $a \in \mathbb{P}$, $\mathbb{N} \subseteq N(n^{-1}a)$. Hence $n\mathbb{N} \subseteq N(a)$, for all $a \in \mathbb{P}$; that is, $n\mathbb{N} \subseteq \mathbb{N}$. (4) is 3.1.3. (5) is immediate from (3).

(3.3) *Let $x^* \in G^* \setminus \{1^*\}$, $\alpha \in \{x^*, C_{x^*}\}$ and set $\mathbb{N} = \mathbb{N}_\alpha$ and $U = U_\alpha$. Then*

(1) $U = \{n \in N : n\mathbb{N} = \mathbb{N}\} = \{n \in N : n\bar{\mathbb{N}} = \bar{\mathbb{N}}\}$.
(2) $U = \{n \in N : \mathbb{N}n = \mathbb{N}\} = \{n \in N : \bar{\mathbb{N}}n = \bar{\mathbb{N}}\}$.
(3) *$U$ is a subgroup of $G$; further, if $\alpha = C_{x^*}$, then $U$ is normal in $G$.*
(4) *If $-1 \in \mathbb{N}$, then $-1 \in U$.*

*Proof.* We start with a proof of (1). Clearly since $N$ is a disjoint union of $\mathbb{N}$ and $\bar{\mathbb{N}}$, $\{n \in N : n\mathbb{N} = \mathbb{N}\} = \{n \in N : n\bar{\mathbb{N}} = \bar{\mathbb{N}}\}$. Let $u \in U$; then by 3.2.3, $u\mathbb{N} \subseteq \mathbb{N}$ and $u^{-1}\mathbb{N} \subseteq \mathbb{N}$. Hence $u\mathbb{N} = \mathbb{N}$. Conversely let $n \in N$ and suppose $n\mathbb{N} = \mathbb{N}$. As $1 \in \mathbb{N}$, $n \in \mathbb{N}$ and as $n^{-1}\mathbb{N} = \mathbb{N}, n^{-1} \in \mathbb{N}$, so $n \in U$. This proves (1). The proof of (2) is identical to the proof of (1). (3) follows from (1) and the fact that if $\alpha = C_{x^*}$, $\mathbb{N}$ is a normal subset of $G$. (4) is immediate from the definition of $U$.

(3.4) *Let $x^* \in G^* \setminus \{1^*\}$ and set $\mathbb{P} = \mathbb{P}_{x^*}$, $U = U_{x^*}$. Let $a \in Nx$ and $n \in N$. Then $n \in N(a)$ if and only if $(nU) \cup (Un) \subseteq N(a)$.*

*Proof.* If $(nU) \cup (Un) \subseteq N(a)$, then since $1 \in U$, $n \in N(a)$. Suppose $n \in N(a)$. Then $1 \in N(n^{-1}a) \cap N(an^{-1})$, by 1.8.1. Hence, by definition, $n^{-1}a$, $an^{-1} \in \mathbb{P}$, so that $U \subseteq N(n^{-1}a) \cap N(an^{-1})$. Now 1.8.1 implies that $(nU) \cup (Un) \subseteq N(a)$, as asserted.

(3.5) *Let $x^* \in G^* \setminus \{1^*\}$ and set $U = U_{x^*}$. Suppose that $U = U_{(x^{-1})^*}$ and that $-1 \in U$. Let $x_1 \in \mathbb{O}_{x^*}$. Then $\mathbb{O}_{x^*} \supseteq (Ux_1) \cup (x_1U)$.*

*Proof.* Let $u \in U$. Suppose $-1 \in N(ux_1)$. Then $-u^{-1} \in N(x_1)$. By 3.4, $U \subseteq N(x_1)$, and in particular, $-1 \in N(x_1)$, a contradiction. Similarly $-1 \notin N(x_1^{-1}u)$, so that $Ux_1 \subseteq \mathbb{O}_{x^*}$. The proof that $x_1U \subseteq \mathbb{O}_{x^*}$ is similar.

(3.6) *Let $x^* \in G^* \setminus \{1^*\}$. Then the following conditions are equivalent.*
(1) $\mathbb{O}_{x^*} = \emptyset$.
(2) *For all $a \in Nx$, $-1 \in N(a) \cup N(a^{-1})$.*
(3) *For all $a \in Nx$, and $n \in N \setminus N(a)$, $a + n \in Nx$.*
(4) *There exists $a \in Nx$ such that for all $n \in N \setminus N(a)$, $a + n \in Nx$.*

*Proof.* (1) if and only if (2) is by definition.

(2) $\to$ (3). Let $a \in Nx$ and $n \in N \setminus N(a)$. Then $-1 \notin N(-n^{-1}a)$; so by (2), $-1 \in N(-a^{-1}n)$; that is, $n^{-1} \in N(a^{-1})$. Hence $a^{-1} + n^{-1} \in N$ and multiplying by $a$ on the right and $n$ on the left we get $a + n \in Na = Nx$.

(3) $\to$ (4). This is immediate.

(4) $\to$ (3). Let $b \in Nx$ and write $b = ma$, for some $m \in N$. Then $N(b) = mN(a)$. Let $n \in N \setminus N(b)$; then $n \notin mN(a)$, so $m^{-1}n \notin N(a)$. Hence, by (4), $a + m^{-1}n \in Nx$, so that $ma + n \in Nx$; that is, $b + n \in Nx$, so (3) holds.

(3) $\to$ (2). Let $a \in Nx$, and suppose $-1 \notin N(a)$. Then by (3), $a - 1 \in Na$. Now, multiplying by $a^{-1}$ on the right we see that $a^{-1} - 1 \in N$; that is, $-1 \in N(a^{-1})$.

(3.7) *Let $a, b \in G \setminus N$ and $\varepsilon \in \{1, -1\}$. Then*

(1) *If $a + b \neq 0$ and $N(a+b) \not\subseteq N(a)$, then*

$$a^*(a+n)^*b^* \text{ is a path in } \Delta, \text{ for any } n \in N(a+b) \setminus N(a).$$

(2) *If $a + b \notin N$ and $N(a) \not\subseteq N(a+b)$, then*

$$b^*(a+b+n)^*(a+b)^* \text{ is a path in } \Delta, \text{ for any } n \in N(a) \setminus N(a+b).$$

(3) *If $a^*z^*(a + \varepsilon b)^*$ is a path in $\Delta$, $\varepsilon \notin N(a^{-1}b)$ and $a^{-1}b \notin N$, then $a^*z^*(\varepsilon + a^{-1}b)^*(a^{-1}b)^*$ is a path in $\Delta$; so in particular, $\mathrm{d}(a^*, (a^{-1}b)^*) \leq 3$.*

*Proof.* For (1), set $c = a+b$ and let $n \in N(a+b) \setminus N(a)$. Then $(a+n)+b = c+n \in N$. By Remark 2.2, $\mathrm{d}((a+n)^*, b^*) \leq 1 \geq \mathrm{d}((a+n)^*, a^*)$, and (1) follows.

For (2), note that $a = (a+b) - b$, so (2) follows from (1).

Finally, for (3), note that $a + \varepsilon b = \varepsilon a(\varepsilon + a^{-1}b)$. Further, $z^*$ commutes with $a^*$ and $(a + \varepsilon b)^*$, so that $z^*$ commutes with $(\varepsilon + a^{-1}b)^*$, and of course $a^{-1}b$ commutes with $(\varepsilon + a^{-1}b)$. Hence, if $(\varepsilon + a^{-1}b)$, $a^{-1}b \notin N$, $a^*z^*(\varepsilon + a^{-1}b)^*(a^{-1}b)^*$ is a path in $\Delta$.

(3.8) *Let $x, y \in G \setminus N$ and let $\bar{n} \in N \setminus (N(x) \cup N(y))$. Suppose $\mathrm{d}(x^*, y^*) > 2$. Then $\bar{n} + m \in N$, for all $m \in N(x + \bar{n}) \cap N(y + \bar{n})$.*

*Proof.* Let $m \in N(x + \bar{n}) \cap N(y + \bar{n})$. Then $x + (\bar{n} + m) \in N$ and $y + (\bar{n} + m) \in N$. Suppose $\bar{n} + m \notin N$. Then, by Remark 2.2, $\mathrm{d}((x^*, (\bar{n}+m)^*) \leq 1 \geq \mathrm{d}(y^*, (\bar{n}+m)^*)$. It follows that $\mathrm{d}(x^*, y^*) \leq 2$, a contradiction.

(3.9) *Let $x^*, y^* \in G^* \setminus \{1^*\}$. Then each of the following conditions imply* $\mathrm{In}(x^*, y^*)$.

(1) $\mathrm{d}(x^*, y^*) > 4$.

(2) $\mathrm{d}(x^*, y^*) > 2$, and $\mathrm{d}(x^*, (x^{-1}y)^*) > 3$.

*Proof.* The fact that (1) implies $\text{In}(x^*, y^*)$ derives from Remark 2.4. Now suppose (2) holds. Let $(a, b) \in Nx \times Ny$. Note that since $\text{d}(a^*, b^*) > 2$, 2.3.1 implies that

(i) $$N(a+b) \subseteq N(a) \cap N(b).$$

Suppose $N(b) \neq N(a+b) \neq N(a)$. Then $N(a) \not\subseteq N(a+b)$ and $N(b) \not\subseteq N(a+b)$, so by 3.7.2,

(ii) $b^*(a+b+n)^*(a+b)^*$ is a path in $\Delta$, for any $n \in N(a) \setminus N(a+b)$

$a^*(a+b+m)^*(a+b)^*$ is a path in $\Delta$, for any $m \in N(b) \setminus N(a+b)$.

From (ii) we get that

(iii) $$a^*(a+b+m)^*(a+b)^*(a+b+n)^*b^* \text{ is a path in } \Delta$$

for any $m \in N(b) \setminus N(a+b)$ and $n \in N(a) \setminus N(a+b)$. Suppose $1 + a^{-1}b \in N$, then $(a+b)^* = a^*$, and then from (iii) we get that $\text{d}(a^*, b^*) \leq 2$, contradicting the choice of $a^*, b^*$. Hence $1 + a^{-1}b \notin N$, so by 3.7.3, $\text{d}(a^*, (a^{-1}b)^*) \leq 3$, a contradiction.

We may now conclude that either $N(a+b) = N(a)$, or $N(a+b) = N(b)$. Hence, by (i), either $N(a) \subseteq N(b)$, or $N(b) \subseteq N(a)$, as asserted.

(3.10) *Let $x^*, y^* \in G^* \setminus \{1^*\}$ and assume $\text{In}(x^*, y^*)$. Then either $\text{Inc}(y^*, x^*)$ or $\text{Inc}(x^*, y^*)$.*

*Proof.* Suppose that $\text{Inc}(y^*, x^*)$ is false. Then, there exists $b \in \mathbb{P}_{y^*}$, such that $N(a) \not\supseteq N(b)$, for all $a \in \mathbb{P}_{x^*}$. Thus $\text{Inc}(x^*, y^*)$ holds.

(3.11) *Let $x^*, y^* \in G^* \setminus \{1^*\}$ such that $\text{In}(x^*, y^*)$. Then*
1. *If $\text{Inc}(y^*, x^*)$, then $\mathbb{N}_{y^*} \supseteq \mathbb{N}_{x^*}$, and $U_{y^*} \geq U_{x^*}$.*
2. *If $(a, b) \in Nx \times Ny$ such that $N(b) \supseteq N(a)$, then $N(-b) \supseteq N(a)$.*
3. *If $(a, b) \in Nx \times Ny$ such that $N(b) \supsetneq N(a)$, then $N(-b) \supsetneq N(a)$.*
4. *If $\text{Inc}(y^*, x^*)$, then $-1 \in \mathbb{N}_{y^*}$ and hence $-1 \in U_{y^*}$.*

*Proof.* For (1), let $b \in \mathbb{P}_{y^*}$. By $\text{Inc}(y^*, x^*)$, there exists $a \in \mathbb{P}_{x^*}$ such that $N(b) \supseteq N(a)$. But, by definition, $N(a) \supseteq \mathbb{N}_{x^*}$. Hence $N(b) \supseteq \mathbb{N}_{x^*}$. As this holds for all $b \in \mathbb{P}_{y^*}$, $\mathbb{N}_{y^*} \supseteq \mathbb{N}_{x^*}$. Then, it is immediate from the definition of $U_{x^*}$ that $U_{y^*} \geq U_{x^*}$.

Let $(a, b) \in Nx \times Ny$ such that $N(b) \supseteq N(a)$. Let $s \in N(b)$. Suppose $-s \notin N(b)$. Then $-s \notin N(a)$ and $-s \in N(-b)$. Hence, by $\text{In}(x^*, y^*)$, $N(-b) \supsetneq N(a)$. Thus we may assume that $-s \in N(b)$, for all $s \in N(b)$. But then $N(-b) = N(b)$, by 1.8.1, and again $N(-b) \supseteq N(a)$; in addition, if $N(b) \supsetneq N(a)$, then $N(-b) = N(b) \supsetneq N(a)$. This show (2) and (3).

Suppose Inc($y^*, x^*$). Let $b \in \mathbb{P}_{y^*}$; then there exists $a \in \mathbb{P}_{x^*}$, such that $N(b) \supseteq N(a)$. By (2), $N(b) \supseteq N(-a)$, so as $-1 \in N(-a)$, $-1 \in N(b)$, as this holds for all $b \in \mathbb{P}_{y^*}$, $-1 \in \mathbb{N}_{y^*}$. This proves the first part of (4) and the second part of (4) is immediate from the definitions.

(3.12) *Let $x^*, y^* \in G^* \setminus \{1^*\}$ and assume*

(i) $\mathrm{d}(x^*, y^*) > 2$.
(ii) $\mathrm{In}(x^*, y^*)$.

*Let $(a, b) \in Nx \times Ny$ and suppose $N(b) \supseteq N(a)$. Then*

(1) $N(a + \varepsilon b) = N(a)$, *for $\varepsilon \in \{1, -1\}$.*
(2) *If $N(b) \supsetneq N(a)$, then $a^*(a + \varepsilon b + n_\varepsilon)^*(a + \varepsilon b)^*$ is a path in $\Delta$, for any $n_\varepsilon \in N(\varepsilon b) \setminus N(a)$, where $\varepsilon \in \{1, -1\}$.*

*Proof.* First note that by 3.11.2, $N(-b) \supseteq N(a)$. Let $\varepsilon \in \{1, -1\}$. As $\mathrm{d}(a^*, b^*) > 2$, $N(a + \varepsilon b) \subseteq N(a)$, by 2.3.1. Let $m \in N(a)$. Then $m \in N(\varepsilon b)$. Suppose $m \notin N(a + \varepsilon b)$. Consider the element $z = a + \varepsilon b + m$. Since $m \notin N(a + \varepsilon b)$, $z \notin N$. However, since $z = a + (\varepsilon b + m)$ (and $\varepsilon b + m \in N$), Remark 2.2 implies that $\mathrm{d}(z^*, a^*) \leq 1$. Similarly as $z = \varepsilon b + (a + m)$ (and $a + m \in N$), $\mathrm{d}(z^*, b^*) \leq 1$. Thus $\mathrm{d}(a^*, b^*) \leq 2$, a contradiction. This shows (1).

Assume $N(b) \supsetneq N(a)$. Then by 3.11.3, $N(-b) \supsetneq N(a)$. Let $n_\varepsilon \in N(\varepsilon b) \setminus N(a)$; then (2) follows from 3.7.2.

(3.13) *Let $x^*, y^* \in G^* \setminus \{1^*\}$ and assume that $\mathrm{d}(x^*, y^*) > 3 < \mathrm{d}(x^*, (x^{-1}y)^*)$. Let $x_1 \in \mathbb{O}_{x^*}$ and $b \in Ny$, such that $1 \notin N(b)$. Then $N(x_1) \supseteq N(b)$.*

*Proof.* First note that by 3.9, $\mathrm{In}(x^*, y^*)$. Suppose $N(x_1) \subsetneq N(b)$. Then, by 3.12, $N(x_1 - b) = N(x_1)$, and is a path in $\Delta$, for any $s \in N(b) \setminus N(x_1)$. Suppose $x_1^{-1}b + 1 \in N$; that is, $1 \in N(x_1^{-1}b)$. Then $-1 \in N(-x_1^{-1}b)$, so $N(x_1^{-1}(-b)) \not\subseteq N(x_1^{-1})$. As $-1 \notin N(-b)$, 2.9.1 implies that $\mathrm{d}(x_1^*, b^*) \leq 3$, contradicting $\mathrm{d}(x^*, y^*) > 3$. Thus $1 \notin N(x_1^{-1}b)$. Hence by 3.7.3, $\mathrm{d}(x^*, (x^{-1}y)^*) \leq 3$, a contradiction.

(*) $$x_1^*(x_1 + b + s)^*(x_1 + b)^*$$

(3.14) *Let $x^*, y^* \in G^* \setminus \{1^*\}$ and assume one of the following conditions holds*

(1) $\mathrm{d}(x^*, y^*) > 4$.
(2) $\mathrm{d}(x^*, y^*) > 3$, $\mathrm{In}(x^*, y^*)$ *and either $\mathbb{O}_{x^*} = \emptyset$ or $\mathbb{O}_{y^*} = \emptyset$.*

*Then, $\mathrm{T}(x^*, y^*)$.*

*Proof.* If $d(x^*, y^*) > 4$, then by 3.9, $\text{In}(x^*, y^*)$ holds. Let $(a, b) \in Nx \times Ny$ and let $\bar{n} \in N \setminus (N(a) \cup N(b))$. By $\text{In}(x^*, y^*)$, we may assume without loss of generality that $N(b) \supseteq N(a)$. By 3.12, $N(a - b) = N(a)$. Note that if (2) holds, then, by 3.6, either $a + \bar{n} \in Na$, or $b + \bar{n} \in Nb$; hence, in any case, by Remark 2.2, $d(a + \bar{n}, b + \bar{n}) > 2$. But $a - b = (a + \bar{n}) - (b + \bar{n})$, and then 2.3.1 implies that $N(a + \bar{n}) \supseteq N(a - b) = N(a)$. Further, by 3.11, $N(-b) \supseteq N(a)$, and as $-\bar{n} \notin N(-b)$, $-\bar{n} \notin N(a)$. Also $a + b = (a - \bar{n}) + (b + \bar{n})$, and if (2) holds, then by 3.6, either $a - \bar{n} \in Na$, or $b + \bar{n} \in Nb$. Hence again, in any case $d(a - \bar{n}, b + \bar{n}) > 2$ and as above we get $N(b + \bar{n}) \supseteq N(a + b) = N(a)$. This shows $\text{T}(x^*, y^*)$.

(3.15) *Let $x^*, y^* \in G^* \setminus \{1^*\}$. Suppose that*
(a) $d(x^*, y^*) > 2$.
(b) $-1 \in \mathbb{N}_{y^*}$.
(c) *For all $\bar{n} \in \bar{\mathbb{N}}_{C_{y^*}}$ and $m \in \mathbb{N}_{C_{x^*}}$, $\bar{n} + m \in N$.*

*Then $G$ satisfies the $U$-Hypothesis with respect to $\mathbb{N}_{C_{y^*}}$.*

*Proof.* Set $\mathbb{N} = \mathbb{N}_{C_{y^*}}$ and $\mathbb{P} = \mathbb{P}_{C_{y^*}}$. First note that by (b) and 3.1.4, $-1 \in \mathbb{N}$. We first claim that

(i) $\qquad\qquad b + m \in \mathbb{N}_{C_{x^*}},\ \text{for all}\ b \in \mathbb{P}\ \text{and}\ m \in \mathbb{N}_{C_{x^*}}.$

To prove (i), let $b \in \mathbb{P}$ and $m \in \mathbb{N}_{C_{x^*}}$. Let $a \in \mathbb{P}_{C_{x^*}}$. Suppose $a + m \in \bar{\mathbb{N}}$. Then, by 3.2.5, $-a - m \in \bar{\mathbb{N}}$, and by (c), $(-a - m) + m \in N$; hence $-a \in N$, a contradiction. Thus $a + m \in \mathbb{N}$ and hence $b + (a + m) \in N$. We have shown that

(ii) $\qquad a + (b + m) \in N,\ \text{for all}\ a \in \mathbb{P}_{C_{x^*}}, b \in \mathbb{P}\ \text{and}\ m \in \mathbb{N}_{C_{x^*}}.$

Since $d(x^*, y^*) > 2$, we can choose $a_1 \in \mathbb{P}_{C_{x^*}}$ so that $d(a_1^*, b^*) > 2$ (see 1.8.5). By (ii), given $m \in \mathbb{N}_{C_{x^*}}$, $a_1 + (b + m) \in N$, so if $b + m \notin N$, then by Remark 2.2, $d(a_1^*, (b + m)^*) \leq 1 \geq d(b^*, (b + m)^*)$, so $d(a_1^*, b^*) \leq 2$, a contradiction. This shows that $b + m \in N$. Now (ii) implies (i). Next we claim:

(iii) $\qquad\qquad\text{For all}\ \bar{n} \in \bar{\mathbb{N}}\ \text{and}\ m \in \mathbb{N}_{C_{x^*}}, \bar{n} + m \in \mathbb{N}.$

Let $b \in \mathbb{P}$, $\bar{n} \in \bar{\mathbb{N}}$ and $m \in \mathbb{N}_{C_{x^*}}$. By (i), $b + m \in \mathbb{N}_{C_{x^*}}$, and by (c), $b + m + \bar{n} \in N$. As this holds for all $b \in \mathbb{P}$, $m + \bar{n} \in \mathbb{N}$, and (iii) is proved.

Finally, let $\bar{n} \in \bar{\mathbb{N}}$. Then by 3.2.5, $-\bar{n} \in \bar{\mathbb{N}}$, and since $1 \in \mathbb{N}_{C_{x^*}}$, (iii) implies that $-\bar{n} + 1 \in \mathbb{N}$. Hence

(iv) $\qquad\qquad\qquad\qquad \bar{n} - 1 \in \mathbb{N}.$

Now (iii), (iv), our assumption (b) and 3.2 imply that $G$ satisfies the $U$-Hypothesis with respect to $\mathbb{N}$.

(3.16) THEOREM. *Let $x^*, y^* \in G^* \setminus \{1^*\}$. Suppose that*

(a) $\mathrm{d}(x^*, y^*) > 2$.
(b) $-1 \in \mathbb{N}_{y^*}$.
(c) *For all $\bar{n} \in \bar{\mathbb{N}}_{y^*}$ and $m \in \mathbb{N}_{x^*}$, $\bar{n} + m \in N$.*

*Then*

(1) *For all $\bar{n} \in \mathbb{N}_{C_{y^*}}$ and $m \in \mathbb{N}_{C_{x^*}}$, $\bar{n} + m \in N$.*
(2) *$G$ satisfies the $U$-Hypothesis with respect to $\mathbb{N}_{C_{y^*}}$.*

*Proof.* Set $\mathbb{N} = \mathbb{N}_{C_{y^*}}$ and let $\bar{n} \in \bar{\mathbb{N}}$ and $m \in \mathbb{N}_{C_{x^*}}$. We want to show that $\bar{n} + m \in N$. After conjugation with some element of $G$, and using 3.1, we may assume that $\bar{n} \in \bar{\mathbb{N}}_{y^*}$. But $m \in \mathbb{N}_{C_{x^*}} \subseteq \mathbb{N}_{x^*}$, so (1) follows from our assumption (c). Then (2) follows from 3.15.

(3.17) THEOREM. *Let $x^*, y^* \in G^* \setminus \{1^*\}$ and assume*

(i) $\mathrm{d}(x^*, y^*) > 2$.
(ii) $\mathrm{Inc}(y^*, x^*)$ *and* $\mathrm{T}(x^*, y^*)$.

*Then $G$ satisfies the $U$-Hypothesis with respect to $\mathbb{N}_{C_{y^*}}$.*

*Proof.* Set $\mathbb{N} = \mathbb{N}_{C_{y^*}}$. We verify assumptions (b) and (c) of Theorem 3.16. Assumption (b) follows from $\mathrm{Inc}(y^*, x^*)$ and 3.11.4.

It remains to verify assumption (c) of Theorem 3.16. Let $\bar{n} \in \bar{\mathbb{N}}_{y^*}$ and let $m \in \mathbb{N}_{x^*}$. By definition, there exists $b \in \mathbb{P}_{y^*}$, such that $\bar{n} \notin N(b)$. Let $a \in \mathbb{P}_{x^*}$, such that $N(b) \supseteq N(a)$ (using $\mathrm{Inc}(y^*, x^*)$). By $\mathrm{T}(x^*, y^*)$, $N(a + \bar{n}) \supseteq N(a) \subseteq N(b + \bar{n})$. In particular, $m \in \mathbb{N}_{x^*} \subseteq N(a) \subseteq N(a + \bar{n}) \cap N(b + \bar{n})$. Since $\mathrm{d}(x^*, y^*) > 2$, 3.8 implies that $\bar{n} + m \in N$, as asserted.

(3.18) THEOREM. *Suppose that $\mathrm{diam}(\Delta) > 4$. Then there exist conjugacy classes $A^*, B^* \subseteq G^* \setminus \{1^*\}$ such that*

(1) *$G$ satisfies the $U$-Hypothesis with respect to $\mathbb{N}_{B^*}$.*
(2) *For all $b \in \mathbb{P}_{B^*}$, there exists $a \in \mathbb{P}_{A^*}$ such that $\mathrm{d}(a^*, b^*) > 4$ and $N(b) \supseteq N(a)$.*

*Proof.* Let $x^*, y^* \in \Delta$ be such that $\mathrm{d}(x^*, y^*) > 4$. By 3.9, $\mathrm{In}(x^*, y^*)$ and by 3.10, we may assume that $\mathrm{Inc}(y^*, x^*)$. Further by 3.14, $\mathrm{T}(x^*, y^*)$. Set $B^* = C_{y^*}$ and $A^* = C_{x^*}$. By Theorem 3.17, (1) holds. Let $b \in \mathbb{P}_{B^*}$. Then there exists $g \in G$, such that $b^g \in \mathbb{P}_{y^*}$ (see 3.1.1). Since $\mathrm{Inc}(y^*, x^*)$, there exists $a \in \mathbb{P}_{x^*}$ such that $N(b^g) \supseteq N(a)$. By 1.8.2, $N(b) \supseteq N(a^{g^{-1}})$. Of course $a^{g^{-1}} \in \mathbb{P}_{A^*}$ and $\mathrm{d}(b^*, (a^{g^{-1}})^*) > 4$, so (2) holds.

## 4. The proof that if $\Delta$ is balanced then $G$ satisfies the $U$-Hypothesis

In this section we continue the notation and definitions of Sections 2 and 3.

*Definitions.* (1) We define a binary relation $\mathfrak{B}$ on $(G^* \setminus \{1^*\}) \times (G^* \setminus \{1^*\})$ as follows. Let $(x^*, y^*) \in (G^* \setminus \{1^*\}) \times (G^* \setminus \{1^*\})$,

$\mathfrak{B}(x^*, y^*)$: The distances $d(x^*, y^*)$, $d(x^*, x^*y^*)$, $d(y^*, x^*y^*)$, $d(x^*, (x^{-1}y)^*)$, $d(y^*, (x^{-1}y)^*)$ are all greater than 3.

(2) We say that $\Delta$ *is balanced* if there exists $x^*, y^* \in G^* \setminus \{1^*\}$ such that $\mathfrak{B}(x^*, y^*)$.

The purpose of this section is to prove the following theorem.

(4.1) THEOREM. *Suppose that $\Delta$ is balanced. Then there exists a conjugacy class $C^* \subseteq G^* \setminus \{1^*\}$ such that*

(1) *$G$ satisfies the U-Hypothesis with respect to $\mathbb{N}_{C^*}$.*
(2) *One of the following holds*:
   (2a) $\mathbb{O}_{x^*} = \emptyset$, *for some $x^* \in G^* \setminus \{1^*\}$.*
   (2b) *For all $m \in \mathbb{M}_{C^*}$, there exists $z^* \in C^*$, such that $m \in N(z_1)$, for all $z_1 \in \mathbb{O}_{z^*}$.*

(4.2) (1) $\mathfrak{B}$ *is symmetric.*
(2) *If $\mathfrak{B}(x^*, y^*)$, then $\mathfrak{B}((x^{-1})^*, y^*)$.*

*Proof.* Suppose $\mathfrak{B}(x^*, y^*)$. We must show that $\mathfrak{B}(y^*, x^*)$. By definition, $d(y^*, x^*) > 3$. Next since $d(y^*, x^*y^*) > 3$, conjugating with $y^*$ we get that $d(y^*, y^*x^*) > 3$. Since $d(x^*, x^*y^*) > 3$, conjugating with $x^*$ we get $d(x^*, y^*x^*) > 3$. Since $d(y^*, (x^{-1}y)^*)$, inverting $(x^{-1}y)^*$, we see that $d(y^*, (y^{-1}x)^*) > 3$, finally since $d(x^*, (x^{-1}y)^*)$, inverting $(x^{-1}y)^*$, we get that $d(x^*, (y^{-1}x)^*) > 3$. Hence $\mathfrak{B}(y^*, x^*)$. The proof of (2) is similar.

*Notation.* From now until the end of Section 4 we fix $x, y \in G \setminus N$ such that $\mathfrak{B}(x^*, y^*)$. We set
$$S := (\{x, x^{-1}\} \times \{y, y^{-1}\}) \cup (\{y, y^{-1}\} \times \{x, x^{-1}\})$$
and
$$\mathbb{O}_S = \mathbb{O}_{x^*} \cup \mathbb{O}_{(x^{-1})^*} \cup \mathbb{O}_{y^*} \cup \mathbb{O}_{(y^{-1})^*}.$$

(4.3) *Let $(g, h) \in S$, then*
(1) $\mathfrak{B}(g^*, h^*)$.
(2) $\mathrm{In}(g^*, h^*)$.

*Proof.* (1) follows from 4.2 and (2) follows from (1) and 3.9.

(4.4) *Suppose $\mathbb{O}_{x^*} = \emptyset$ or $\mathbb{O}_{y^*} = \emptyset$. Then $G$ satisfies the U-Hypothesis.*

*Proof.* First note that by $\mathfrak{B}(x^*, y^*)$, 4.3 and 3.14, $\mathrm{T}(x^*, y^*)$. Then, by 4.3, and 3.10, we may assume without loss that $\mathrm{Inc}(y^*, x^*)$. Now the lemma follows from Theorem 3.17.

In view of 4.4, and symmetry, we assume from now on that

The sets $\mathbb{O}_{x^*}, \mathbb{O}_{(x^{-1})^*}, \mathbb{O}_{y^*}$, and $\mathbb{O}_{(y^{-1})^*}$ are not empty.

*Notation.* Given $z \in \{x, x^{-1}, y, y^{-1}\}$, $z_1$ will always denote an element in $\mathbb{O}_{z^*}$.

(4.5) *Let* $g \in \{x, x^{-1}, y, y^{-1}\}$; *then* $-1 \in \mathbb{N}_{g^*}$, $-1 \in \mathbb{N}_{C_{g^*}}$ *and* $-1 \in U_{g^*}$.

*Proof.* Let $g \neq h \in \{x, x^{-1}, y, y^{-1}\}$, with $h \notin \{g, g^{-1}\}$. By 4.3.1, $\mathfrak{B}(g^*, h^*)$. It suffices to show that $-1 \in \mathbb{N}_{g^*}$, then by 3.1.4, $-1 \in \mathbb{N}_{C_{g^*}}$, and by 3.3.4, $-1 \in U_{g^*}$. Letting $a \in \mathbb{P}_{g^*}$, we must show that $-1 \in N(a)$. Suppose $-1 \notin N(a)$, then, $1 \notin N(-a)$, so by 3.13, $N(g_1) \supseteq N(-a)$. But $-1 \in N(-a)$, a contradiction.

(4.6) *Let* $z \in \mathbb{O}_S$. *Then*
(1) $N(z) = N(h)$, *for all* $h \in \mathbb{O}_S$.
(2) $1 \notin N(z)$.
(3) *If* $a \in Nz$ *such that* $1 \notin N(a)$, *then* $N(z) \supseteq N(a)$.
(4) *If* $\bar{n} \in \bar{\mathbb{N}}_{z^*}$, *then* $\bar{n}^{-1} \in N(z)$.
(5) $N(z) = \mathbb{M}_{z^*}$.
(6) $\mathbb{N}_{z^*}$ *is independent of the choice of* $z$.
(7) $U_{z^*}$ *is independent of the choice of* $z$.

*Proof.* We show that $\mathfrak{B}(x^*, y^*)$ implies $N(x_1) \supseteq N(y_1)$. Then, (1) follows from 4.3.1. A similar application of 4.3.1 will be used throughout the proof. Now $1 \notin N(-y_1)$, so by 3.13, $N(x_1) \supseteq N(-y_1)$. Then, by 3.11, $N(x_1) \supseteq N(y_1)$.

Suppose $1 \in N(x_1)$. Then $-1 \in N(-x_1)$, so that $N(-x_1) \supsetneq N(y_1)$. By 3.11, $N(x_1) \supsetneq N(y_1)$, contradicting (1). Hence (2) holds.

(3) is immediate from 3.11, (1) and 4.3.2. To show (4), let $\bar{n} \in \bar{\mathbb{N}}_{z^*}$. By definition, there exists $a \in \mathbb{P}_{z^*}$, such that $\bar{n} \notin N(a)$. Then, $1 \notin N(\bar{n}^{-1}a)$, and so by (3), $N(z) \supseteq N(\bar{n}^{-1}a)$. But $\bar{n}^{-1} \in N(\bar{n}^{-1}a)$, so that $\bar{n}^{-1} \in N(z)$.

Next let $h \in \mathbb{O}_S$, with $h^* \neq z^*, (z^{-1})^*$. Note that $N(h) \subseteq N(b)$, for all $b \in \mathbb{P}_{z^*}$, by $\mathrm{In}(h^*, z^*)$, so $N(z) = N(h) \subseteq \mathbb{N}_{z^*}$. Let $u \in U_{z^*}$. If $u \in N(z)$; then, by 3.4, $U_{z^*} \subseteq N(z)$, a contradiction, as $-1 \in U_{z^*}$. Hence $N(z) \subseteq \mathbb{M}_{z^*}$. Let $m \in \mathbb{M}_{z^*}$; then, by definition, $m^{-1} \in \bar{\mathbb{N}}_{z^*}$, so by (4), $m = (m^{-1})^{-1} \in N(z)$. Hence $N(z) = \mathbb{M}_{z^*}$, and (5) holds.

To show (6), by 4.3.1, it suffices to show that $\bar{\mathbb{N}}_{x^*} \subseteq \bar{\mathbb{N}}_{y^*}$ (so $\mathbb{N}_{x^*} \supseteq \mathbb{N}_{y^*}$). Let $\bar{n} \in \bar{\mathbb{N}}_{x^*}$; then by (4) and (1), $\bar{n}^{-1} \in N(y_1)$. But by (5), $N(y_1) = \mathbb{M}_{y^*}$, so by definition, $\bar{n} = (\bar{n}^{-1})^{-1} \in \bar{\mathbb{N}}_{y^*}$. Finally (7) is immediate from (6).

(4.7) *Let $\bar{n} \in \bar{\mathbb{N}}_{x^*}$ and $m \in \mathbb{N}_{y^*}$. Then $\bar{n} + m \in N$.*

*Proof.* Set $\mathbb{N} = \mathbb{N}_{x^*}$, $\mathbb{M} = \mathbb{M}_{x^*}$ and $U = U_{x^*}$. Note that by 4.6, $\mathbb{N} = \mathbb{N}_{z^*}$, $\mathbb{M} = \mathbb{M}_{z^*} = N(z)$, and $U = U_{z^*}$, for all $z \in \mathbb{O}_S$. First we claim that

(i) $$z + \bar{n} \in Nz, \text{ for all } z \in \mathbb{O}_S.$$

Indeed, by 4.6.4, $\bar{n}^{-1} \in \mathbb{M}$, so as $\mathbb{M} = N(z^{-1})$, $z^{-1} + \bar{n}^{-1} \in N$ and (i) holds.

Further, by 3.12, $N(x_1 - y_1) = N(x_1) = \mathbb{M}$, and by (i), $d((x_1 + \bar{n})^*, (y_1 + \bar{n})^*) > 3$, hence, by 2.3.1, $\mathbb{M} = N(x_1 - y_1) = N((x_1 + \bar{n}) - (y_1 + \bar{n})) \subseteq N(x_1 + \bar{n})$. Similarly, $\mathbb{M} \subseteq N(y_1 + \bar{n})$, so that

(ii) $$N(x_1 + \bar{n}) \supseteq \mathbb{M} \subseteq N(y_1 + \bar{n})$$

by 3.8, $\bar{n} + \mathbb{M} \subseteq N$, for all $\bar{n} \in \bar{\mathbb{N}}$. We have shown

(iii) $$\bar{n} + m \in N, \text{ for all } \bar{n} \in \bar{\mathbb{N}} \text{ and } m \in \mathbb{M}.$$

Next we show that $\bar{n} + 1 \in N$, for all $\bar{n} \in \bar{\mathbb{N}}$. We first claim that

(iv) $$N(x_1 + 1) \supseteq \mathbb{M}.$$

Suppose not and let $m \in \mathbb{M} \setminus N(x_1 + 1)$; recall that by 3.12, $N(x_1 - y_1) = N(x_1) = \mathbb{M}$. But $x_1 - y_1 = (x_1 + 1) - (y_1 + 1)$, so $m \in N(x_1 - y_1) \setminus N(x_1 + 1)$. Hence, by 3.7.1,

(v) $$(x_1 + 1)^*(x_1 + 1 + m)^*(y_1 + 1)^* \text{ is a path in } \Delta.$$

Replacing $y_1$, by $y_1^{-1}$, the same argument shows that

(vi) $$(x_1 + 1)^*(x_1 + 1 + m)^*(y_1^{-1} + 1)^* \text{ is a path in } \Delta.$$

It follows from (v) and (vi) that $(x_1 + 1 + m)^*$ commutes with $(y_1^{-1} + 1)$ and $(y_1 + 1)$. But $y_1 + 1 = y_1(y_1^{-1} + 1)$, so $(x_1 + 1 + m)^*$ commutes with $y_1^*$. However, applying Remark 2.2 twice, we see that $d((x_1 + 1 + m)^*, x_1^*) \leq 2$. Hence we get that $d(x_1^*, y_1^*) \leq 3$, contradicting $\mathfrak{B}(x^*, y^*)$. This shows (iv). Similarly, $N(y_1 + 1) \supseteq \mathbb{M}$. Since $\bar{n}^{-1} \in \mathbb{M}$, 3.8 implies that $\bar{n}^{-1} + 1 \in N$, so $\bar{n} + 1 = \bar{n}(\bar{n}^{-1} + 1) \in N$. We have shown

(vii) $$\bar{n} + 1 \in N, \text{ for all } \bar{n} \in \bar{\mathbb{N}}.$$

Let $u \in U$. Then $u^{-1}\bar{n} \in \bar{\mathbb{N}}$, by 3.3, so by (vii), $u^{-1}\bar{n} + 1 \in N$, so $\bar{n} + u \in N$. We have shown

(viii) $$\bar{n} + u \in N, \text{ for all } u \in U.$$

Since $\mathbb{N}$ is the union of $\mathbb{M}$ and $U$, (iii) and (viii) complete the proof.

(4.8) *$G$ satisfies the $U$-Hypothesis with respect to $\mathbb{N}_{C_{x^*}}$.*

*Proof.* This follows immediately from 4.5, 4.7 and Theorem 3.16.

(4.9) *Let* $\mathbb{N} = \mathbb{N}_{C_{x^*}}$ *and* $\mathbb{M} = \mathbb{M}_{C_{x^*}}$. *Then* $\mathbb{N} = \mathbb{N}_{C_{z^*}}$ *and* $\mathbb{M} = \mathbb{M}_{C_{z^*}}$, *for all* $z \in \mathbb{O}_S$.

*Proof.* Let $z \in \mathbb{O}_S$. By definition, $\mathbb{N}_{C_{x^*}} = \bigcap\{\mathbb{N}_{v^*} : v^* \in C_{x^*}\}$ and $\mathbb{N}_{C_{z^*}} = \bigcap\{\mathbb{N}_{v^*} : v^* \in C_{z^*}\}$. But, by 4.6.6 and 3.1.2, $\{\mathbb{N}_{v^*} : v^* \in C_{x^*}\} = \{\mathbb{N}_{v^*} : v^* \in C_{z^*}\}$, so $\mathbb{N} = \mathbb{N}_{C_{z^*}}$. Then, by definition, $\mathbb{M} = \mathbb{M}_{C_{z^*}}$.

(4.10) *Set* $\mathbb{M} = \mathbb{M}_{C_{x^*}}$, *and let* $m \in \mathbb{M}$. *Then there exists* $z^* \in C_{x^*}$, *such that* $m \in N(z_1)$, *for all* $z_1 \in \mathbb{O}_{z^*}$.

*Proof.* Since $m \in \mathbb{M}$, $m \in \mathbb{N}_{C_{x^*}}$. Since $m \notin U_{C_{x^*}}$, there exists $z^* \in C_{x^*}$, such that $m \notin U_{z^*}$. Hence $m \in \mathbb{M}_{z^*}$. After conjugation, and using 3.1, we may assume that $z = x$. But then the lemma follows from 4.6.

Note now that by 4.4, 4.9 and 4.10, Theorem 4.1 holds.

## 5. The $U$-Hypothesis

In this section $\emptyset \neq \mathbb{N} \subsetneq N$ is a proper subset of $N$ such that $\mathbb{N}$ is a normal subset of $G$. We denote $\bar{\mathbb{N}} = N \setminus \mathbb{N}$ and assume the $U$-Hypothesis.

($U1$) $1, -1 \in \mathbb{N}$.
($U2$) $\mathbb{N}^2 = \mathbb{N}$.
($U3$) For all $\bar{n} \in \bar{\mathbb{N}}$, $\bar{n} + 1 \in \mathbb{N}$ and $\bar{n} - 1 \in N$.

(5.1) *Remark.* Notice that if $\operatorname{diam}(\Delta) > 4$ or $\Delta$ is balanced, then by Theorems 3.18 and 4.1, $G$ satisfies the $U$-Hypothesis with respect to $\mathbb{N} = \mathbb{N}_{X^*}$, where $X^* = B^*$, if $\operatorname{diam}(\Delta) > 4$ ($B^*$ as in Theorem 3.18) and $X^* = C^*$ if $\Delta$ is balanced ($C^*$ as in Theorem 4.1).

(5.2) *Let* $U = \{n \in \mathbb{N} : n^{-1} \in \mathbb{N}\}$. *Then*
(1) $U = \{n \in N : n\mathbb{N} = \mathbb{N}\} = \{n \in N : n\bar{\mathbb{N}} = \bar{\mathbb{N}}\}$.
(2) $U = \{n \in N : \mathbb{N}n = \mathbb{N}\} = \{n \in N : \bar{\mathbb{N}}n = \bar{\mathbb{N}}\}$.
(3) $U$ *is a normal subgroup of* $G$.
(4) $-1 \in U$.

*Proof.* This was already proved in 3.3 in a slightly different context; for completeness we include a proof. Clearly since $N$ is a disjoint union of $\mathbb{N}$ and $\bar{\mathbb{N}}$, $\{n \in N : n\mathbb{N} = \mathbb{N}\} = \{n \in N : n\bar{\mathbb{N}} = \bar{\mathbb{N}}\}$. Let $u \in U$, then by ($U2$), $u\mathbb{N} \subseteq \mathbb{N}$ and $u^{-1}\mathbb{N} \subseteq \mathbb{N}$. Hence $u\mathbb{N} = \mathbb{N}$. Conversely let $n \in N$ and suppose $n\mathbb{N} = \mathbb{N}$. As $1 \in \mathbb{N}$, $n \in \mathbb{N}$ and as $n^{-1}\mathbb{N} = \mathbb{N}, n^{-1} \in \mathbb{N}$, so that $n \in U$. This proves (1). The proof of (2) is identical to the proof of (1). (3) follows from (1) and the fact that $\mathbb{N}$ is a normal subset of $G$. Note that (4) follows immediately from ($U1$).

(5.3) *Notation.* We denote $\mathbb{M} = \mathbb{N} \setminus U$. Hence $N = \mathbb{M} \dot{\cup} U \dot{\cup} \bar{\mathbb{N}}$ is a disjoint union.

(5.4) (1) *For all $\bar{n} \in \bar{\mathbb{N}}$, $\bar{n} + U = U$.*

(2) *For all $\bar{n} \in \bar{\mathbb{N}}$, $\bar{n}^{-1} \in \mathbb{N}$.*

*Proof.* We first show

(i) $\qquad\qquad$ For all $\bar{n} \in \bar{\mathbb{N}}, \bar{n} - 1 \in \mathbb{N}$.

Let $\bar{n} \in \bar{\mathbb{N}}$ and suppose $\bar{n} - 1 \notin \mathbb{N}$, then, $\bar{n} - 1 \in \bar{\mathbb{N}}$ and by $(U3)$, $(\bar{n}-1)+1 \in \mathbb{N}$, a contradiction. This shows (i).

Let $\bar{m} \in \bar{\mathbb{N}}$. Suppose that $\bar{m}^{-1} \in \mathbb{N}$, then by $(U2)$, $\bar{m}^{-1}\mathbb{N} \subseteq \mathbb{N}$. We conclude that $\bar{m}^{-1}(\bar{m} \pm 1) \in \mathbb{N}$. Hence $\bar{m}^{-1} \pm 1 \in \mathbb{N}$. Suppose $\bar{m}^{-1} \in \bar{\mathbb{N}}$. Then by $(U3)$ and (i), $\bar{m}^{-1} \pm 1 \in \mathbb{N}$. Hence in either case we get that

(ii) $\qquad\qquad \bar{m}^{-1} \pm 1 \in \mathbb{N}, \text{ for all } \bar{m} \in \bar{\mathbb{N}}.$

Next we show

(iii) $\qquad\qquad$ For all $\bar{n} \in \bar{\mathbb{N}}, \quad \bar{n} \pm 1 \in U.$

Let $\bar{n} \in \bar{\mathbb{N}}$ and let $\varepsilon \in \{1, -1\}$. By (i) and $(U3)$, $\bar{n} + \varepsilon \in \mathbb{N}$. Hence we must show that $(\bar{n} + \varepsilon)^{-1} \in \mathbb{N}$. Suppose $(\bar{n} + \varepsilon)^{-1} \notin \mathbb{N}$. Set $\bar{m} = (\bar{n} + \varepsilon)^{-1}$. Then $\bar{m} \in \bar{\mathbb{N}}$, so by (ii), $\bar{m}^{-1} - \varepsilon \in \mathbb{N}$. But $\bar{m}^{-1} - \varepsilon = \bar{n} \in \bar{\mathbb{N}}$, a contradiction.

We can now prove (1). Let $u \in U$ and $\bar{n} \in \bar{\mathbb{N}}$. Then by 5.2.1, $u^{-1}\bar{n} \in \bar{\mathbb{N}}$ and by (iii), $u^{-1}\bar{n} + 1 \in U$. It follows that $\bar{n} + u = u(u^{-1}\bar{n} + 1) \in U$. Hence

(iv) $\qquad\qquad \bar{n} + U \subseteq U.$

Next by 5.2.4, $-u \in U$, and by (iv), $\bar{n} - u \in U$. Again by 5.2.4, $u - \bar{n} \in U$ and hence $u = \bar{n} + (u - \bar{n}) \in \bar{n} + U$. Hence $U \subseteq \bar{n} + U$ and (1) is proved.

Finally we prove (2). Let $\bar{n} \in \bar{\mathbb{N}}$ and suppose $\bar{n}^{-1} \notin \mathbb{N}$. Then $\bar{n}^{-1} \in \bar{\mathbb{N}}$, so by (1), $\bar{n}^{-1} + 1 \in U$. Then by 5.2.2, $\bar{n} + 1 = \bar{n}(\bar{n}^{-1} + 1) \in \bar{\mathbb{N}}$, which contradicts $(U3)$.

(5.5) (1) *For all $s \in N \setminus U$, $s \in \mathbb{M}$ if and only if $s^{-1} \in \bar{\mathbb{N}}$.*

(2) *For all $\bar{n} \in \bar{\mathbb{N}}$, $\bar{n} + U = U$.*

(3) *For all $u \in U, u\bar{\mathbb{N}} = \bar{\mathbb{N}}u = \bar{\mathbb{N}}$ and $u\mathbb{M} = \mathbb{M}u = \mathbb{M}$.*

(4) $\bar{\mathbb{N}}^2 \subseteq \bar{\mathbb{N}}$ *and* $\mathbb{M}^2 \subseteq \mathbb{M}$.

*Proof.* For (1) let $m \in \mathbb{M} \subseteq \mathbb{N}$. If $m^{-1} \in \mathbb{N}$, then, by definition, $m \in U$, a contradiction. Hence $m^{-1} \in \bar{\mathbb{N}}$. Let $n \in \bar{\mathbb{N}}$. By 5.4.2, $n^{-1} \in \mathbb{N}$, and since $n \notin U, n^{-1} \in \mathbb{M}$. This shows (1). (2) is from 5.4.1 and (3) is from 5.2.1 and 5.2.2.

Let $\bar{n}, \bar{m} \in \bar{\mathbb{N}}$ and suppose $\bar{n}\bar{m} \in \mathbb{N}$. By (1), $\bar{n}^{-1} \in \mathbb{N}$, and by ($U2$), $\bar{m} = \bar{n}^{-1}(\bar{n}\bar{m}) \in \mathbb{N}$, a contradiction. Hence $\bar{\mathbb{N}}^2 \subseteq \bar{\mathbb{N}}$. Let $m, m' \in \mathbb{M}$. Suppose $mm' \in U \cup \bar{\mathbb{N}}$. Then $m^{-1} \in \bar{\mathbb{N}}$ (by (1)) and by (3) and the fact that $\bar{\mathbb{N}}^2 \subseteq \bar{\mathbb{N}}$, $m' = m^{-1}(mm') \in \bar{\mathbb{N}}$, a contradiction. Hence $\mathbb{M}^2 \subseteq \mathbb{M}$.

## 6. Further consequences of the $U$-Hypothesis

In this section we continue the notation and hypotheses of Section 5, deriving further consequences. We denote $\Gamma = N/U$ (note that by 5.2.3, $U$ is a normal subgroup of $G$ and hence of $N$). Recall from 1.3 that we denote by $\nu: G \to F^\#$ the reduced norm function, in the case when $[D:F] < \infty$.

(6.1) *Definition.* We define an order relation $\leq$ on $\Gamma$ as follows. For $Ua, Ub \in \Gamma$, $Ua < Ub$ if and only if $Ua \neq Ub$ and $ba^{-1} \in \bar{\mathbb{N}}$.

(6.2) (1) *The relation $\leq$ is a well defined linear order relation on $\Gamma$.*
(2) *If $Ua, Ub, Uc, Ud \in \Gamma$, with $Ua \leq Uc$ and $Ub \leq Ud$, then $Uab \leq Ucd$.*

*Proof.* It is clear from 5.5.3 that $\leq$ is independent on coset representatives and hence it is a well defined relation on $\Gamma$. We show it is an order relation. If $Ua < Ub$, then $ba^{-1} \in \bar{\mathbb{N}}$; hence by 5.5.1, $ab^{-1} \in \mathbb{M}$ and it follows that $Ub \not< Ua$. Also if $Ua < Ub < Uc$, then $ba^{-1} \in \bar{\mathbb{N}}$ and $cb^{-1} \in \bar{\mathbb{N}}$. Hence by 5.5.4, $ca^{-1} = (cb^{-1})(ba^{-1}) \in \bar{\mathbb{N}}$ and hence $Ua < Uc$. Finally let $Ua, Ub \in \Gamma$, with $Ua \neq Ub$. Then by 5.5.1 either $ab^{-1} \in \bar{\mathbb{N}}$ or $ba^{-1} \in \bar{\mathbb{N}}$; hence either $Ua < Ub$, or $Ub < Ua$, so $\leq$ is linear.

For (2), if $Ua = Uc$, or $Ub = Ud$, then (2) follows directly from the definition of $\leq$ and the fact that $\bar{\mathbb{N}}$ is a normal subset of $G$. So suppose $Ua < Uc$ and $Ub < Ud$. Then $ca^{-1}, db^{-1} \in \bar{\mathbb{N}}$. Now $(cd)(ab)^{-1} = cdb^{-1}a^{-1} = ca^{-1}adb^{-1}a^{-1}$. Since $\bar{\mathbb{N}}$ is a normal subset of $G$, $adb^{-1}a^{-1} \in \bar{\mathbb{N}}$. By 5.5.4, $\bar{\mathbb{N}}^2 \subseteq \bar{\mathbb{N}}$, so $ca^{-1}adb^{-1}a^{-1} \in \bar{\mathbb{N}}$. Hence, $(cd)(ab)^{-1} \in \bar{\mathbb{N}}$ and $Uab < Ucd$, as asserted.

(6.3) *Let $Ua, Ub \in \Gamma$, with $Ua \neq Ub$. Then*
(1) $Ua + Ub \subseteq N$, *and*
(2) $Ua + Ub = \min\{Ua, Ub\}$.

*Proof.* Without loss of generality we may assume that $Ua < Ub$. Let $x \in Ua$ and $y \in Ub$. Then $yx^{-1} \in \bar{\mathbb{N}}$. Hence by 5.5.2, $1 + yx^{-1} \in U$ and multiplying by $x$ on the right we see that $x + y \in Ux = Ua$. This shows (1) and the fact that $Ua + Ub \subseteq Ua$. But $Ua + Ub$ contains the coset $U(a+b)$, and it follows that $Ua + Ub = Ua$.

(6.4) COROLLARY. *Let $Ua \in \Gamma$ and let $x, y \in Ua$. Suppose $x + y \in N$. Then $U(x+y) \geq Ua$.*

*Proof.* Suppose $U(x+y) < Ua = Ux$. Then by 6.3, $y = (x+y) - x \in U(x+y)$. But $y \in Ua$, a contradiction.

(6.5) COROLLARY. *Let $a_1, a_2, \ldots, a_k \in N$ and assume there exists some $1 \leq i \leq k$, such that $Ua_i < Ua_j$, for all $j \neq i$. Then $Ua_1 + Ua_2 + \cdots + Ua_k = Ua_i$.*

*Proof.* This follows immediately from 6.3 by induction.

(6.6) *Suppose $[D:F] < \infty$ and let $n \in N \setminus UF^\#$. Then there exists $r \leq \deg(D)$ such that $n^r \in UF^\#$.*

*Proof.* Let
$$\alpha_0 + \alpha_1 x^{k_1} + \cdots + \alpha_t x^{k_t}$$
be the minimal polynomial of $n$ over $F$ with $\alpha_i \neq 0$, for all $0 \leq i \leq t$ and $0 < k_1 < k_2 < \cdots < k_t$. Suppose there exists some $0 \leq i \leq t$, such that $U\alpha_i n^{k_i} < U\alpha_j n^{k_j}$, for all $j \neq i$. Then by 6.5, $\alpha_0 + \alpha_1 n^{k_1} + \cdots + \alpha_t n^{k_t} \in U\alpha_i n^{k_i}$. In particular, $\alpha_0 + \alpha_1 n^{k_1} + \cdots + \alpha_t n^{k_t} \neq 0$, a contradiction. Hence the set of minimal elements in the set $\{U\alpha_0, U\alpha_1 n^{k_1}, \ldots, U\alpha_t n^{k_t}\}$ is of size larger than 1. It follows that there are indices $0 \leq i < j \leq t$, such that $U\alpha_i n^{k_i} = U\alpha_j n^{k_j}$. We conclude that $n^{k_j - k_i} \in U(\alpha_i \alpha_j^{-1})$. Note now that $r = k_j - k_i \leq k_t \leq \deg(D)$ and that $n^r \in UF^\#$.

(6.7) *Suppose $[D:F] < \infty$ and let $n \in N$, with $\nu(n) \in U$. Then $n \in U$.*

*Proof.* Suppose first that $n \in UF^\#$. Note that as $U \triangleleft G$, for each $u \in U$, $\nu(u) \in U$. This is because $\nu(u)$ is a product of conjugates of $u$ (see 1.4). Write $n = \alpha u$, with $u \in U$ and $\alpha \in F^\#$. Then
$$\nu(n) = \alpha^{\deg(D)} \nu(u)$$
and it follows that $\alpha^{\deg(D)} = \nu(n)\nu(u)^{-1} \in U$. By 5.5.4, $\alpha \in U$, and hence $n \in U$.

Next suppose $n \in N \setminus UF^\#$. Then by (6.6), $n^r \in UF^\#$, for some $1 < r \leq \deg(D)$. Note now that $\nu(n^r) = \nu(n)^r \in U$, so by the previous paragraph of the proof, $n^r \in U$, this contradicts 5.5.4.

(6.8) COROLLARY. *If $[D:F] < \infty$, then $N/U \leq Z(G/U)$.*

*Proof.* Here $Z(G/U)$ is the center of $G/U$. Let $g \in G$ and $n \in N$. Then $\nu([g,n]) = 1 \in U$. Hence by 6.7, $[g,n] \in U$.

(6.9) *Remark.* Note that if $[D : F] < \infty$, then the canonical homomorphism $v : N \to \Gamma$ behaves like a *valuation* on $N$ in the sense that $v$ is a group homomorphism, and $\Gamma$ is a linearly ordered abelian group. Further $v(a + b) \geq \min\{v(a), v(b)\}$, whenever $a + b \in N$. In particular the restriction $v : F^\# \to v(F^\#)$ is a valuation on $F$.

(6.10) *If $[D : F] < \infty$, then $F^\# \not\subseteq U$.*

*Proof.* Suppose $F^\# \subseteq U$ and let $\bar{n} \in \bar{\mathbb{N}}$. Let

$$\alpha_0 + \alpha_1 x^{k_1} + \cdots + \alpha_t x^{k_t}$$

be the minimal polynomial of $\bar{n}$ over $F$ with $\alpha_i \neq 0$, for all $0 \leq i \leq t$ and $0 < k_1 < k_2 < \cdots < k_t$. Then

$$\alpha_0 + \alpha_1 \bar{n}^{k_1} + \cdots + \alpha_t \bar{n}^{k_t} = 0.$$

We show by induction on $j$ that $\alpha_0 + \alpha_1 \bar{n}^{k_1} + \cdots + \alpha_j \bar{n}^{k_j} \in U$, for all $0 \leq j \leq t$. By hypothesis $\alpha_0 \in U$. Suppose $\alpha_0 + \alpha_1 \bar{n}^{k_1} + \cdots + \alpha_j \bar{n}^{k_j} \in U$. Note that as $\alpha_{j+1} \in U$, 5.5.3 and 5.5.4 imply that $\alpha_{j+1} \bar{n}^{k_{j+1}} \in \bar{\mathbb{N}}$; hence by 5.5.2, $(\alpha_0 + \alpha_1 \bar{n}^{k_1} + \cdots + \alpha_j \bar{n}^{k_j}) + \alpha_{j+1} \bar{n}^{k_{j+1}} \in U$. But we cannot have $\alpha_0 + \alpha_1 \bar{n}^{k_1} + \cdots + \alpha_t \bar{n}^{k_t} \in U$, a contradiction.

## 7. Towards the proof of Theorem A

In this and the following sections we finally prove Theorem A. We continue the notation of the previous sections. In particular, $\Delta$ is the commuting graph of $G^*$. We assume that either $\text{diam}(\Delta) > 4$, or $\Delta$ is balanced. If $\text{diam}(\Delta) > 4$ then we fix $A^*, B^*$ to denote the conjugacy classes as in Theorem 3.18. Recall that $\hat{A} = \{a \in G : a^* \in A^*\}$ and $\hat{B} = \{b \in G : b^* \in B^*\}$. If $\Delta$ is balanced, then we fix $C^*$ to denote the conjugacy class as in Theorem 4.1; again $\hat{C} = \{c \in G : c^* \in C^*\}$.

If $\text{diam}(\Delta) > 4$, let $X^* = B^*$, while if $\Delta$ is balanced let $X^* = C^*$. We let $\mathbb{P} = \mathbb{P}_{X^*}$, $\mathbb{N} = \mathbb{N}_{X^*}, \bar{\mathbb{N}} = \bar{\mathbb{N}}_{X^*}$, $\mathbb{M} = \mathbb{M}_{X^*}$ and $U = U_{X^*}$. Note that by Remark 5.1, all the results of Sections 5 and 6 apply here.

In this section we further assume that $G^*$ is a nonabelian finite simple group and that $[D : F] < \infty$. We draw the attention of the reader to Remarks 2.2 and 2.4.

*Definitions and Notation.* (1) $\hat{K} = \{a \in \mathbb{O}U \setminus N : N(a) \supseteq \mathbb{M}\}$.
(2) $K^* = \{a^* : a \in \hat{K}\}$.
(3) An element $a \in G \setminus N$ is a *standard element* if it satisfies the following condition: If $n \in N(a)$, then $Un \subseteq N(a)$.
(4) We denote by $\Phi$ the set of all standard elements in $G \setminus N$.

(7.1) (1) $G = \mathbb{O}N$.
(2) $(\mathbb{O}U) \cap N = U$.
(3) $\mathbb{O}U/U \simeq G^*$.
(4) $[G, N] \leq U$.

*Proof.* (1) follows from our assumption that $G^*$ is simple and from 1.5. Let $n \in (\mathbb{O}U) \cap N$; then $n = au$, for some $a \in \mathbb{O}$, so $\nu(n) = \nu(u)$. Since $U$ is normal in $G$, $\nu(u) \in U$, by 1.4. Then by 6.7, $n \in U$. Next, since $G = (\mathbb{O}U)N, G^* = G/N \simeq \mathbb{O}U/(\mathbb{O}U) \cap N = \mathbb{O}U/U$, by (2). Finally, (4) is from 6.8.

(7.2) *Let $a, b \in G \setminus N$. Then*
(1) Let $n \in N(a)$, then $a + n \in Un$.
(2) Let $n \in N(a)$, then $Um \subseteq N(a)$, for all $Um < Un$. Further if $a \in \Phi$, then also $Un \subseteq N(a)$.
(3) Let $n \in N \setminus N(a)$, then $Um \leq Un$, for all $m \in N(a)$. Further if $a \in \Phi$, then $Um < Un$, for all $m \in N(a)$.
(4) If $a \in \Phi$ and $b \in G \setminus N$, then $N(a) \subseteq N(b)$ or $N(b) \subseteq N(a)$.
(5) Let $n \in N$. Then $\mathbb{N} \subseteq N(n)$ if and only if $n \in \bar{\mathbb{N}}$ and $\mathbb{M} \subseteq N(n)$ if and only if $n \in U \cup \bar{\mathbb{N}}$.

*Proof.* For (1), suppose $a + n = m \notin Un$. Note that as $-1 \in U$, $-n \in Un$ and hence $a = m - n \in Um + Un \subseteq N$, by 6.3.1, a contradiction.

For (2), assume $Um < Un$. By (1), $a = n + nu$, for some $u \in U$. Then $a + m = n + nu + m \in Um$, by 6.5. Hence $m \in N(a)$. This proves the first part of (2) and the second part of (2) is obvious. Now (3) is an immediate consequence of (2).

Let $a \in \Phi$ and $b \in G \setminus N$ and suppose $N(b) \not\subseteq N(a)$. Let $n \in N(b) \setminus N(a)$; then by (2), $Um \subseteq N(b)$, for all $Um < Un$. By (3), if $m \in N(a)$, then $Um < Un$. Hence, $N(a) \subseteq N(b)$. This proves (4).

We now prove (5). Let $n \in \mathbb{N}$. Then $-n \notin N(n)$, by definition, and $-n \in \mathbb{N}$, thus $\mathbb{N} \not\subseteq N(n)$. Let $\bar{n} \in \bar{\mathbb{N}}$; then by 6.3, $\mathbb{N} \subseteq N(n)$. This proves the first part of (5). The proof of the second part of (5) is similar.

(7.3) *Let $a \in G \setminus N$. Then*
(1) If $a \in \Phi$, then $Na \subseteq \Phi$.
(2) $a \in \Phi$ if and only if

(∗)          For each $b \in Na$ such that $1 \in N(b), U \subseteq N(b)$.

In particular, if $a \notin \Phi$, then there exists $b \in Na$, with $N(b) \cap U \neq \emptyset$, but $U \not\subseteq N(b)$.
(3) If $a \in \Phi$, then $N(a)$ is a normal subset of $G$.
(4) $\Phi$ is a normal subset of $G$.

*Proof.* Suppose $a \in \Phi$. Let $b \in Na$ and $m \in N(b)$. Write $b = sa$, $s \in N$. Then $s^{-1}m \in N(a)$. Let $u \in U$. Since $a \in \Phi$, $s^{-1}mu \in N(a)$, so that $mu \in N(sa) = N(b)$. Thus $Um \subseteq N(b)$ as asserted.

For (2), note that (1) implies that if $a \in \Phi$, then $(*)$ holds. So assume $(*)$ holds. Let $m \in N(a)$. Then $1 \in N(m^{-1}a)$, so by $(*)$, $U \subseteq N(m^{-1}a)$. Hence $Um \subseteq N(a)$ and $a \in \Phi$.

Next let $a \in \Phi$, $n \in N(a)$ and $g \in G$. By 7.1.4, $n^g \in Un \subseteq N(a)$, so $N(a)$ is a normal subset of $G$. (4) follows from (3) since for $a \in \Phi$ and $g \in G$, $N(a^g) = g^{-1}N(a)g = N(a)$; so $a^g \in \Phi$.

(7.4) (1) *If* $\mathrm{diam}(\Delta) > 4$, *then* $\hat{B} \subseteq \Phi$.

(2) *If* $\Delta$ *is balanced, then* $\hat{C} \subseteq \Phi$.

*Proof.* (1) and (2) follow from the definition of $U$ and from 7.3.2.

(7.5) *Let* $a \in \Phi$. *Then*

(1) *For all* $u \in U, N(ua) = N(a) = N(au)$.

(2) *If* $n \in N \setminus N(a)$, *then* $N(a+n) \supseteq N(a)$.

(3) *Let* $x, y \in (\mathbb{O}U) \cap Na$. *Then* $N(x) = N(y)$.

*Proof.* For (1) note that $N(ua) = uN(a) \subseteq N(a)$, as $a \in \Phi$. Similarly $u^{-1}N(a) \subseteq N(a)$, so $N(a) \subseteq uN(a) = N(ua)$. This proves the first part of (1) and the proof of the second part of (1) is the same. For (2), let $m \in N(a)$; then, by 7.2.3, $Um < Un$, and by 6.3.2, $a + n + m = a + um$, for some $u \in U$. Then, since $a \in \Phi$, $a + um \in N$, so that $m \in N(a+n)$.

Next we prove (3): notice that $xy^{-1} \in (\mathbb{O}U) \cap N$. Hence, by 7.1.2, $xy^{-1} \in U$, so (3) follows from (1).

(7.6) *Let* $a \in G \setminus N$ *and* $m \in N$. *Suppose* $Um \subseteq N(x)$, *for some* $x \in \hat{C}_{a^*} \cap (\mathbb{O}U)$; *then* $Um \subseteq N(z)$, *for all* $z \in \hat{C}_{a^*} \cap (\mathbb{O}U)$.

*Proof.* Recall that $C_{a^*}$ is the conjugacy class of $a^*$ in $G^*$ and $\hat{C}_{a^*} = \{c \in G : c^* \in C_{a^*}\}$. First we claim that

$(*)$ $\qquad\qquad\qquad Um \subseteq N(x^g)$, for all $g \in G$.

Let $g \in G$. Then, by 1.8.2, $N(x^g) = g^{-1}N(x)g$. Hence $N(x^g) \supseteq g^{-1}(Um)g = Um^g = Um$, where the last equality follows from 7.1.4.

Let $z \in \hat{C}_{a^*} \cap (\mathbb{O}U)$. Then there exists $g \in G$, such that $z^* = (x^g)^*$. By $(*)$, $Um \subseteq N(x^g)$. Hence, we may assume that $Nx = Nz$. But then $xz^{-1} \in N \cap (\mathbb{O}U) = U$ (see 7.1.2). Hence, there exists $u \in U$ such that $z = ux$. Then $N(z) = uN(x)$, so that $N(z) \supseteq u(Um) = Um$, as asserted.

(7.7) *Let $a, b \in G \setminus N$. Then*

(1) *If $U \subseteq N(a) \cap N(b)$, then $U \subseteq N(ab)$.*
(2) *If $U \cap N(a) \neq \emptyset$ and $\mathbb{M} \subseteq N(b)$, then $\mathbb{M} \subseteq N(ab)$.*
(3) *If $U \subseteq N(a) \cap N(b)$, then $U \subseteq N(a+b)$.*
(4) *If $\mathbb{M} \subseteq N(a) \cap N(b)$, then $\mathbb{M} \subseteq N(a+b)$.*
(5) *Suppose $a \in \mathbb{O}U \setminus N$ and let $\ell > 1$, such that $a^\ell \in N$. Then $a^\ell \in U$.*
(6) *Suppose $a \in \mathbb{O}U \setminus N$. Then $U \not\subseteq N(a)$. In particular, if $a \in \Phi$, then $N(a) \subseteq \mathbb{M}$.*

*Proof.* For (1), let $u \in U$. Then $ab + u = ab - b + b + u = (a-1)b + (b+u)$. As $-1 \in U$, $a - 1 \in U$, by 7.2.1. Further by 7.2.1, $b + u \in U$, write $v = a - 1$ and $w = b + u$. Then $ab + u = vb + w = v(b + v^{-1}w) \in N$. Hence $u \in N(ab)$.

For (2), let $u \in U \cap N(a)$ and let $m \in \mathbb{M}$. Then $ab + m = ab + ub + (m - ub) = (a + u)b + (m - ub)$. Note now that by 7.2.1, $a + u = v$ and $m - ub = wm$, for some $v, w \in U$. Hence $ab + m = vb + wm = v(b + v^{-1}wm) \in N$, where the last equality is because $\mathbb{M} \subseteq N(b)$ and because $(v^{-1}w)\mathbb{M} = \mathbb{M}$.

For (3), let $u \in U$; then $(a + b) + u = a + (b + u)$. But since $u \in N(b)$, $b + u = v \in U$, by 7.2.1. Hence $(a + b) + u = a + v \in N$. Thus $U \subseteq N(a+b)$. The proof of (4) is similar.

Assume the hypotheses of (5). Since $a \in \mathbb{O}U$, $\nu(a) \in U$, so $\nu(a^\ell) \in U$. Hence by 6.7, $a^\ell \in U$. Let $a \in \mathbb{O}U \setminus N$. Since $G^*$ is finite there exists $r \geq 2$, with $a^r \in N$. By (5), $a^r \in U$. Hence $a^{-1} = ua^{r-1}$, for some $u \in U$. Suppose $U \subseteq N(a)$. Then by (1), $U \subseteq N(a^{r-1})$, so $U \subseteq N(a^{-1})$. In particular $1 \in N(a) \cap N(a^{-1})$ contradicting 1.8.4. The second part of (6) follows from the first part of (6) and by 7.2.2.

(7.8) *Let $s \in \mathbb{M}$ and suppose that*

(∗) $\qquad\qquad\qquad s^2 \in N(z)$, *for all* $z \in \mathbb{O}U \setminus N$.

*Then $s \in N(z)$, for all $z \in \mathbb{O}U \setminus N$.*

*Proof.* Assume that there exists $x \in \mathbb{O}U \setminus N$, such that $s \notin N(x)$. Set $y := -s^{-1}x$. Then $-1 \notin N(y)$. First we claim that

(∗∗) $\qquad\qquad\qquad -1 \in N(yy^g)$, *for all* $g \in G$.

This is because $yy^g = (s^{-1}x)(s^{-1})^g x^g = s^{-2}xx^g u$, for some $u \in U$, where the last equality follows from 7.1.4. Since $x \in \mathbb{O}U$, $-xx^g u \in \mathbb{O}U$, so, if $-xx^g u \notin N$, then, by hypothesis (∗), $s^2 \in N(-xx^g u)$. If $-xx^g u \in N$, then $-xx^g u \in (\mathbb{O}U) \cap N = U$ (see 7.1.2). Since $s \in \mathbb{M}, s^2 \in \mathbb{M}$, by 5.5.4, so, by 7.2.5, $s^2 \in N(-xx^g u)$ in this case too. Now in any case $s^2 \in N(-xx^g u)$, and it follows that $-1 \in N(yy^g)$.

Now, taking $a = y = b$ in 2.10, we get from 2.10 and $(**)$, that $G^*$ is not simple, a contradiction.

## 8. Some properties of $\hat{K}$ and the proof that $\hat{K} \neq \emptyset$

In this section we continue the notation and hypotheses of Section 7 recalling from there that we defined

$$\hat{K} = \{a \in \mathbb{O}U \setminus N : N(a) \supseteq \mathbb{M}\}.$$

(8.1) (1) $\hat{K}$ *is a normal subset of* $G$.
(2) *If* $a \in \hat{K}$, *then* $Ua \subseteq \hat{K}$.

*Proof.* (1) follows immediately from the fact that $\mathbb{M}$ and $\mathbb{O}U$ are normal subsets of $G$ and from 1.8.2. (2) follows from the fact that $u\mathbb{M} = \mathbb{M}$, from 1.8.1 and the definition of $\hat{K}$.

(8.2) *Suppose there exists* $a \in \mathbb{O}U \setminus N$ *such that* $N(a) \cap U \neq \emptyset$. *Then*
(1) *For all* $b \in \hat{K}$ *such that* $ab \in G \setminus N$, $ab \in \hat{K}$.
(2) $\hat{K} = \mathbb{O}U \setminus N$.

*Proof.* First note that by 7.2.2, $\mathbb{M} \subseteq N(a)$, so that $a \in \hat{K}$. Next, for (1), let $b \in \hat{K}$. Then $N(b) \supseteq \mathbb{M}$. By 7.7.2, $\mathbb{M} \subseteq N(ab)$, and clearly $ab \in \mathbb{O}U$, hence $ab \in \hat{K}$.

Next, since $\hat{K}$ is a normal subset of $G$, $\hat{K}^* \cup \{1^*\}$ is a normal subset of $G^*$. Further note that by (1), $a^*(\hat{K}^* \cup \{1^*\}) \subseteq \hat{K}^* \cup \{1^*\}$. Hence, by 1.9, $\hat{K}^* \cup \{1^*\} = G^*$. Let $b \in \mathbb{O}U \setminus N$. Then $b^* = k^*$, for some $k \in \hat{K}$, and then $bk^{-1} \in (\mathbb{O}U) \cap N \leq U$. Hence $b = uk$, for some $u \in U$, so $b \in \hat{K}$. It follows that $\hat{K} = \mathbb{O}U \setminus N$.

(8.3) *Assume that* $\mathrm{diam}(\Delta) > 4$ *and that for all* $m \in \mathbb{M}$, *there exists* $z \in (\hat{A} \cup \hat{B}) \cap (\mathbb{O}U)$ *such that* $Um \subseteq N(z)$. *Then* $\hat{K} \neq \emptyset$.

*Proof.* Let $\mathbb{V} = \bigcap_{x \in \hat{A} \cap \mathbb{O}U} N(x)$ and $\mathbb{W} = \bigcap_{y \in \hat{B} \cap \mathbb{O}U} N(y)$. Let $m \in \mathbb{V}$, $u \in U$ and $x \in \hat{A} \cap (\mathbb{O}U)$. Then $u^{-1}x \in \hat{A} \cap (\mathbb{O}U)$, so $m \in N(u^{-1}x)$. Thus $um \in N(x)$ and $Um \subseteq N(x)$. As this holds for all $x \in \hat{A} \cap (\mathbb{O}U)$, $Um \subseteq \mathbb{V}$. Similarly, $Um \subseteq \mathbb{W}$, for all $m \in \mathbb{W}$. Next, if $Um \subseteq \mathbb{V}$ and $Us \leq Um$, for some $s \in N$, then, by 7.2.2, $Us \subseteq \mathbb{V}$. Similarly if $Um \subseteq \mathbb{W}$ and $Us \leq Um$, for some $s \in N$, then $Us \subseteq \mathbb{W}$.

Next we claim that either $\mathbb{V} \subseteq \mathbb{W}$, or $\mathbb{W} \subseteq \mathbb{V}$. Suppose $\mathbb{V} \not\subseteq \mathbb{W}$. Let $Um \subseteq \mathbb{V}$, such that $Um \cap \mathbb{W} = \emptyset$. Then, by the previous paragraph of the proof, $Us < Um$, for all $Us \subseteq \mathbb{W}$. Hence, by the previous paragraph of the proof, $Us \subseteq \mathbb{V}$ and hence $\mathbb{W} \subseteq \mathbb{V}$.

Finally, by 7.6, and by the hypothesis of the lemma, $\mathbb{M} \subseteq \mathbb{V} \cup \mathbb{W}$, so, by the second paragraph of the proof $\mathbb{M} \subseteq \mathbb{V}$, or $\mathbb{M} \subseteq \mathbb{W}$. Hence either $\hat{A} \cap (\mathbb{O}U) \subseteq \hat{K}$, or $\hat{B} \cap (\mathbb{O}U) \subseteq \hat{K}$ and $\hat{K} \neq \emptyset$.

(8.4) THEOREM. $\hat{K} \neq \emptyset$.

*Proof.* Suppose $\hat{K} = \emptyset$. By 8.2, we may assume

(∗) $\qquad\qquad U \cap N(x) = \emptyset,$ for all $x \in \mathbb{O}U \setminus N.$

*Case* 1. $\mathrm{diam}(\Delta) > 4$. We shall show that for all $m \in \mathbb{M}$, there exists $z \in (\hat{A} \cup \hat{B}) \cap (\mathbb{O}U)$ such that $Um \subseteq N(z)$. Then, by 8.3, $\hat{K} \neq \emptyset$, a contradiction. Let $m \in \mathbb{M}$. Since $m^{-1} \in \bar{\mathbb{N}}$, there exists $b \in \mathbb{P}$, such that $m^{-1} \notin N(b)$. By 3.18.2, there exists $a \in \mathbb{P}_{A^*}$ such that $N(a) \subseteq N(b)$ and $\mathrm{d}(a^*, b^*) > 4$. Note that by 3.9, $\mathrm{In}(a^*, b^*)$. Further, since $b \in \Phi$, $-m^{-1} \notin N(b)$ (see 7.2.2), and hence $-m^{-1} \notin N(a)$. Let $x \in Na \cap (\mathbb{O}U)$ and $y \in Nb \cap (\mathbb{O}U)$ and suppose that $Um \not\subseteq N(x)$ and $Um \not\subseteq N(y)$. Since $y \in \Phi$, $m \notin N(y)$ and, after replacing $x$ by $ux$, for some $u \in U$, we may assume that $m \notin N(x)$.

Suppose first that $N(y) \supseteq N(x)$. Let $a' = ma$. Notice that $m \in N(a')$ and $-1 \notin N(a')$. Write $a' = xn$, $n \in N$. Notice that $mn^{-1} \in N(x)$, so $mn^{-1} \in N(y)$. Thus $m \in N(yn)$. But $y \in \Phi$, so $n^{-1}N(y)n = N(y)$ (see 7.3.3); thus $m \in N(ny)$. Note now that by (∗), all the hypotheses of 2.11 are met, for $x, y, m$ and $n$, so by 2.11, $\mathrm{d}(x^*, y^*) \leq 4$, contradicting $\mathrm{d}(a^*, b^*) > 4$.

Suppose next that $N(x) \supseteq N(y)$. Let $b' = mb$. Notice that $m \in N(b')$ and $-1 \notin N(b')$. Write $b' = ny$, $n \in N$. Notice that $n^{-1}m \in N(y)$ and since $y \in \Phi$, $mn^{-1} \in N(y)$. Thus $mn^{-1} \in N(x)$ and hence, $m \in N(xn)$. Again we see that by (∗), all the hypotheses of 2.11 are met, for $x, y, m$ and $n$; so by 2.11, $\mathrm{d}(x^*, y^*) \leq 4$, a contradiction. Hence, either $Um \subseteq N(x)$, or $Um \subseteq N(y)$. This completes the proof of the theorem, in the case when $\mathrm{diam}(\Delta) > 4$.

*Case* 2. $\Delta$ is balanced. We use Theorem 4.1. First note that by (∗) and 3.6.2, we are in case (2b) of Theorem 4.1. Let $m \in \mathbb{M}$. By Theorem 4.1, there exists $z \in \hat{C}$ such that $m \in N(z_1)$, for all $z_1 \in \mathbb{O}_{z^*}$. By (∗), $Nz \cap (\mathbb{O}U) \subseteq \mathbb{O}_{z^*}$; thus $m \in N(z_1)$, for some $z_1 \in Nz \cap (\mathbb{O}U)$. Since $z_1 \in \Phi$, $Um \subseteq N(z_1)$ and hence, by 7.6, $Um \subseteq N(x)$, for all $x \in \hat{C} \cap (\mathbb{O}U)$. As this holds for all $m \in \mathbb{M}$, $\hat{C} \cap (\mathbb{O}U) \subseteq \hat{K}$ and $\hat{K} \neq \emptyset$.

## 9. The proof that $\hat{K} = \mathbb{O}U \setminus N$

In this section we continue the notation and hypotheses of Section 7. Note that by Theorem 8.4, $\hat{K} \neq \emptyset$. The purpose of this section is to prove

(9.1) THEOREM. $\hat{K} = \mathbb{O}U \setminus N$.

In view of 8.2, we may (and do) assume that $N(a) \cap U = \emptyset$, for all $a \in \mathbb{O}U \setminus N$.

(9.2) *Suppose that for all $m \in \mathbb{M}$ there exists $s \in \mathbb{M}$, with $Um < Us$. Then*
(1) *Let $a_1, b_1 \in \hat{K}$ such that $a_1 b_1 \in G \setminus N$. Then $a_1 b_1 \in \hat{K}$.*
(2) $\hat{K} = \mathbb{O}U \setminus N$.

*Proof.* For (1), let $m \in \mathbb{M}$ and let $s \in \mathbb{M}$, with $Um < Us$. Then
$$a_1 b_1 + m = a_1 b_1 - a_1 s + (a_1 s + m) = a_1(b_1 - s) + (a_1 s + m).$$
By 7.2.1, $b_1 - s = us$, for some $u \in U$. Next $a_1 s + m = (a_1 + ms^{-1})s$. Note that as $Um < Us$, $ms^{-1} \in \mathbb{M}$, and hence $a_1 + ms^{-1} \in N$. Hence $a_1 s + m \in N$, so by 7.2.1, $a_1 s + m = vm$, for some $v \in U$. Hence we get that $a_1 b_1 + m = a_1(us) + vm$ and as above $a_1(us) + vm \in N$, so $m \in N(a_1 b_1)$. Hence $N(a_1 b_1) = \mathbb{M}$. Since $a_1 b_1 \in \mathbb{O}U \setminus N$, $a_1 b_1 \in \hat{K}$.

The proof of (2) is exactly like the proof of 8.2.2; all we need is the property established in (1).

*Notation.* We fix the letter $m$ to denote an element $m \in \mathbb{M}$ such that $Us \leq Um$, for all $s \in \mathbb{M}$ (see 9.2).

(9.3) (1) *Let $k, \ell \in \mathbb{Z}$ such that $0 < k \leq \ell$. Suppose $x, y \in \mathbb{O}U \setminus N$ such that $N(x) \supseteq Um^k$ and $N(y) \supseteq Um^\ell$. Then $N(xy) \supseteq Um^{k+\ell}$.*
(2) *There exists $t > 0$, such that for all $z \in \mathbb{O}U \setminus N$, $N(z) \supseteq Um^t$.*

*Proof.* For (1), we have
$$xy + m^{k+\ell} = xy + xm^\ell - xm^\ell + m^{k+\ell} = uxm^\ell - xm^\ell + m^{k+\ell}$$
$$= (ux - x + m^k)m^\ell = (ux + vm^k)m^\ell \in N.$$
Here, $u, v \in U$ and we used 7.2.1 for the equalities.

For (2), let $x \in \hat{K}$. Let $X^*$ be the conjugacy class of $x^*$ in $G^*$. Let $\hat{X} = \{x \in G \setminus N : x^* \in X^*\}$. Note that by 7.6, $\hat{X} \cap \mathbb{O}U \subseteq \hat{K}$. Now $G^* = \langle X^* \rangle$, and every element $g^* \in G^*$ can be written as a product of elements in $X^*$. For $g^* \in G^*$, let $\ell(g^*)$ be the minimal length of a word in the alphabet $X^*$ which equals $g^*$. Let $t = \max\{\ell(g^*) : g^* \in G^*\}$. Note that every element in $\mathbb{O}U \setminus N$ can be written as a word of length at most $t$ in the alphabet $\hat{X} \cap (\mathbb{O}U)$. Hence, by (1), as $Um \subseteq N(y)$, for $y \in \hat{X} \cap \mathbb{O}U$, $Um^t \subseteq N(z)$, for all $z \in \mathbb{O}U \setminus N$.

We now complete the proof of Theorem 9.1. Suppose $\hat{K} \neq \mathbb{O}U \setminus N$. Let $1 \leq t \in \mathbb{Z}$, minimal subject to $Um^t \subseteq N(z)$, for all $z \in \mathbb{O}U \setminus N$. Since $\hat{K} \neq \mathbb{O}U \setminus N$, $t \geq 2$. Since there exists $y \in \mathbb{O}U \setminus N$ such that $Um^{t-1} \not\subseteq N(y)$, we may assume without loss of generality that $m^{t-1} \notin N(y)$. Set $s = m^{t-1}$. Notice that $s^2 = m^{2(t-1)}$, and as $2(t-1) \geq t$, we conclude that $s^2 \in N(z)$, for all $z \in \mathbb{O}U \setminus N$. But now, by 7.8, $s \in N(z)$, for all $z \in \mathbb{O}U \setminus N$. This implies that $Us \subseteq N(z)$, for all $z \in \mathbb{O}U \setminus N$, a contradiction.

## 10. The construction of the local ring $R$ and the proof of Theorem A

In this section we continue the hypotheses of Section 7. In addition, in view of Theorem 9.1, we know that $\hat{K} = \mathbb{O}U \setminus N$. We will construct a local ring $R$ and finally prove Theorem A.

(10.1) *Let $a \in G$. Then*
(1) *If $a \notin N$, then $\mathbb{M} \subseteq N(a)$ if and only if $a = na_1$, for some $n \in U \cup \bar{\mathbb{N}}$ and some $a_1 \in \hat{K}$.*
(2) *If $a \notin N$, then $U \subseteq N(a)$ if and only if $a = \bar{n}a_1$, for some $\bar{n} \in \bar{\mathbb{N}}$ and some $a_1 \in \hat{K}$.*
(3) *If $a \in N$, then $\mathbb{M} \subseteq N(a)$, if and only if $a \in U \cup \bar{\mathbb{N}}$.*
(4) *If $a \in N$, then $\mathbb{N} \subseteq N(a)$ if and only if $a \in \bar{\mathbb{N}}$.*

*Proof.* Note first that if $a \notin N$, then by Theorem 9.1, and by 7.1.1, $a = na_1$, for some $n \in N$ and some $a_1 \in \hat{K}$.

Suppose $a \notin N$. Write $a = na_1$, with $n \in N$ and $a_1 \in \hat{K}$. Now suppose that $\mathbb{M} \subseteq N(a)$ and let $u \in U$ such that $u \notin N(a_1)$ (see 7.7.6). Then $a + nu = n(a_1 + u) \notin N$. Hence $nu \notin \mathbb{M}$, so $nu \in U \cup \bar{\mathbb{N}}$. It follows that $n \in U \cup \bar{\mathbb{N}}$. Suppose that $U \subseteq N(a)$; then $Un^{-1} \subseteq N(a_1)$. But by 7.7.6, $U \not\subseteq N(a_1)$ and hence, by 7.2.2, $n^{-1} \in \mathbb{M}$, so that $n \in \bar{\mathbb{N}}$.

Conversely, let $a_1 \in \hat{K}$ and $n \in U \cup \bar{\mathbb{N}}$. If $n \in U$, then $na_1 \in \hat{K}$, so that $\mathbb{M} \subseteq N(na_1)$. If $n \in \bar{\mathbb{N}}$, then for all $u \in U$, $na_1 + u = n(a_1 + n^{-1}u)$, and as $n^{-1}u \in \mathbb{M}$, $na_1 + u \in N$. Hence $U \subseteq N(na_1)$. This completes the proof of (1) and (2). (3) and (4) are as in 7.2.5.

*Definition.* We define
$$R = \{x \in D : x = 0, \text{ or } \mathbb{M} \subseteq N(x)\},$$
$$I = \{r \in R : r = 0 \text{ or } U \subseteq N(r)\}.$$

(10.2) (1) $R \cap N = U \cup \bar{\mathbb{N}}$.

(2) $R \cap (G \setminus N) = \{nk : n \in U \cup \bar{\mathbb{N}} \text{ and } k \in \hat{K}\}$.
(3) $R$ is a subring of $D$.
(4) $I$ is the unique maximal ideal of $R$.
(5) $R \setminus I = \mathbb{O}U$ is the group of unites of $R$.

*Proof.* (1) and (2) are as in 10.1.3 and 10.1.1 respectively. Let $x, y \in R^{\#}$. To show $x + y \in R$, suppose $x \neq -y$. Assume first that $x, y \in N$. If $x + y \in N$, then by 6.3, 6.4 and (1), $x + y \in U \cup \bar{\mathbb{N}}$, so $x + y \in R$. Suppose $x + y \notin N$. Then since $-x \in N(x+y)$, and $-x \in U \cup \bar{\mathbb{N}}$, we get from 7.2.2 that $\mathbb{M} \subseteq N(x+y)$, so $x + y \in R$.

Now assume $x \notin N$. If $y \in N(x)$, then $x + y \in Uy$, by 7.2.1, and as $y \in U \cup \bar{\mathbb{N}}$, $Uy \subseteq U \cup \bar{\mathbb{N}}$, so $x + y \in R$. If $y \in N \setminus N(x)$, then since $y \in U \cup \bar{\mathbb{N}}$, $y + m \in \mathbb{M}$, for all $m \in \mathbb{M}$ and hence $x + y + m \in N$, for all $m \in \mathbb{M}$. Hence $\mathbb{M} \subseteq N(x+y)$, so $x + y \in R$.

Suppose $x, y \notin N$; then by 7.7.4, $x + y \in R$. Let $x, y \in R^{\#}$. It is easy to see that $xy \in R$ by (1) and (2).

The proof of (4) is similar to the proof of (3) from (1), (2), 10.1 and 7.7.3, and we omit the details. Let $r \in R \setminus I$. Then $\mathbb{M} \subseteq N(r) \not\supseteq U$, so by 10.1.1 and 10.1.2, if $r \notin N$, then $r = uk$, for some $k \in \hat{K}$ and $u \in U$; so $r \in \mathbb{O}U$, while if $r \in N$, then by 10.1.3 and 10.1.4, $r \in U$ which shows that $R \setminus I \subseteq \mathbb{O}U$. The inclusion $\mathbb{O}U \subseteq R \setminus I$ follows from the fact that $\mathbb{O}U \subseteq R$ and from 7.7.6. This proves (5).

Let
$$\phi : R \to R/I$$
be the canonical homomorphism. Let
$$\psi : \mathbb{O}U \to (R/I)^{\#}$$
be the multiplicative group homomorphism induced by $\phi$.

(10.3) (1) $R/I$ is a division algebra.
(2) $\psi$ is surjective and $\ker \psi \leq U$.
(3) $R/I$ is infinite.

*Proof.* (1) and the first part of (2) are obvious. Let $r \in \ker \psi$. Then $r - 1 \in I$. Hence $r - 1 = a \in I$. But then $r = a + 1$, and as $a \in I$, $a + 1 \in N$; thus $r \in N$. It follows that $r \in (\mathbb{O}U) \cap N = U$, by 7.1.2. Next we prove (3). Since $\ker \psi \leq U$ and since, by 7.1.3, $\mathbb{O}U/U \simeq G^*$, we see that $G^*$ is a homomorphic image of $\mathbb{O}U/\ker\psi$ ($\mathbb{O}U/U \simeq (\mathbb{O}U/\ker\psi)/(U/\ker\psi)$). Hence $G^*$ is a homomorphic image of $(R/I)^{\#}$. But if $R/I$ is finite, then $R/I$ is a field, which is impossible. Hence $R/I$ is infinite.

(10.4) $\mathbb{O}U \subseteq U + U$.

*Proof.* We apply Theorem 1.6 to the division ring $R/I$. Since $\psi(U)$ is a subgroup of finite index in $(R/I)^{\#}$, Theorem 1.6 implies that for all $r \in R \setminus I$, there are $u_1, u_2 \in U$ such that $r + I = u_1 - u_2 + I$. Hence $r = u_1 - u_2 + a$, with $a \in I$. Note now that by 6.3.2, 7.2.1, 10.1.4 and the definition of $I$, $-u_2 + a \in U$. Hence for all $r \in R \setminus I$, $r = u + v$, with $u, v \in U$. But by 10.2.5, $R \setminus I = \mathbb{O}U$; so the proof is complete.

We can now reach the final contradiction and complete the proof of Theorem A. Note that by 7.3.1 and 7.4, $\hat{K} \cap \Phi \neq \emptyset$ and by 7.7.6, if $k \in \hat{K} \cap \Phi$, then $N(k) = \mathbb{M}$. Let $k \in \hat{K} \cap \Phi$. By 10.4, there are $u, v \in U$, with $k = u + v$. Thus, $-u \in N(k)$. But $N(k) = \mathbb{M}$, a contradiction.

BEN-GURION UNIVERSITY BEER-SHEVA 84105, ISRAEL
*E-mail address*: yoavs@math.bgu.ac.il